\documentclass[10pt]{article}
\usepackage{setspace}
\usepackage{amsmath} 
\usepackage{amssymb}
\usepackage{latexsym}
\usepackage{xcolor}
\usepackage{fancyhdr}
\allowdisplaybreaks 

 \usepackage{comment}

\def\XXint#1#2#3{{\setbox0=\hbox{$#1{#2#3}{\int}$} 
\vcenter{\hbox{$#2#3$}}\kern-.5\wd0}}   

 \numberwithin{equation}{section}
\newtheorem{theorem}[equation]{Theorem}
\newtheorem{proposition}[equation]{Proposition}
\newtheorem{definition}[equation]{Definition}
\newtheorem{remark}[equation]{Remark}
\newtheorem{lemma}[equation]{Lemma}
\newtheorem{corollary}[equation]{Corollary}

\title{
A nonvariational form of the acoustic single layer potential
 } 
 
\author{  
Massimo Lanza de Cristoforis
\\
Dipartimento di Matematica `Tullio Levi-Civita', 
\\
Universit\`a degli Studi di Padova, 
\\
Via Trieste 63, Padova 35121, 
Italy. 
\\
E-mail: mldc@math.unipd.it   }

\date{\ }

\begin{document}
 
 \maketitle

\noindent
{\bf Abstract:}  We consider a bounded open subset $\Omega$ of ${\mathbb{R}}^n$ of class $C^{1,\alpha}$ for some 
$\alpha\in]0,1[$ and the space $V^{-1,\alpha}(\partial\Omega)$ of (distributional) normal derivatives on the boundary of $\alpha$-H\"{o}lder continuous functions in $\Omega$ that have Laplace operator in the  Schauder space with negative exponent  $C^{-1,\alpha}(\overline{\Omega})$. Then we  
 prove those properties of the acoustic single layer potential that are necessary to analyze the Neumann problem for the Helmholtz equation in $\Omega$ with boundary data in $V^{-1,\alpha}(\partial\Omega)$ and solutions in the space 
 of $\alpha$-H\"{o}lder continuous functions in $\Omega$ that have Laplace operator 
 in $C^{-1,\alpha}(\overline{\Omega})$, \textit{i.e.}, in a space of functions
  that may have infinite Dirichlet integral. Namely, a Neumann problem 
 that  does not belong to the classical variational setting. 

 \vspace{\baselineskip}



\section{Introduction} Our starting point are the classical examples of Prym \cite{Pr1871} and  Hadamard \cite{Ha1906} of harmonic functions in the unit ball of the plane that are continuous up to the boundary and have  infinite Dirichlet integral, \textit{i.e.},  whose
 gradient is not square-summable. Such functions solve the classical Dirichlet problem in the unit ball, but not the corresponding weak (variational) formulation. For a discussion on this point we refer to 
Maz’ya and Shaposnikova \cite{MaSh98}, Bottazzini and Gray \cite{BoGr13} and
 Bramati,  Dalla Riva and  Luczak~\cite{BrDaLu23} which contains examples of H\"{o}lder continuous harmonic functions with infinite Dirichlet integral.

In the papers \cite{La24b},\cite{La24c}, the author has analyzed the Neumann problem for the Laplace and the Poisson equation and in \cite{La25} the author has considered the uniqueness problem for $\alpha$-H\"{o}lder continuous solutions of exterior boundary value problems for the Helmoltz equation that satisfy the Sommerfeld radiation condition and that  may have infinite Dirichlet integral close to the boundary. Then 
the present paper is a further step in order to prove the existence for the
  the Neumann problem for the Helmholtz equation. Since we plan to exploit the layer potential theoretic method, we plan to analyze the properties of the single layer potential that are necessary to investigate the Neumann problem. Here the main difference from the classical and variational treatment  is that the single layers we consider may well have infinite Dirichlet integrals around the boundary. 

We consider a bounded open subset $\Omega$ of ${\mathbb{R}}^n$ of class $C^{1,\alpha}$, the space $C^{-1,\alpha}(\overline{\Omega})$ of sums of $\alpha$-H\"{o}lder continuous functions and of first order partial distributional derivatives of $\alpha$-H\"{o}lder continuous functions in $\Omega$     and 
we exploit a distributional normal derivative $\partial_\nu$ on $\partial\Omega$  for functions $u$ in the space $C^{0,\alpha}(\overline{\Omega})_\Delta$ of the $\alpha$-H\"{o}lder continuous functions of $C^{0,\alpha}(\overline{\Omega})$ such that the distributional Laplace operator $\Delta u$ belongs to $C^{-1,\alpha}(\overline{\Omega})$ (as in \cite[\S 5]{La24c}).  If $u\in  C^{0,\alpha}(\overline{\Omega})_\Delta$, then $\partial_\nu u$ belongs to a subspace
\[
V^{-1,\alpha}(\partial\Omega) 
\]
of the dual of $C^{1,\alpha}(\partial\Omega)$ (cf.~\cite[Thm.~5.10]{La24c}). In order to show the solvability of the Neumann problem for the Helmholtz equation in $\Omega$ with solutions in $C^{0,\alpha}(\overline{\Omega})_\Delta$ and
  Neumann datum $g\in V^{-1,\alpha}(\partial\Omega)$,  	one would have to show  the membership of an acoustic  single layer in the space $C^{0,\alpha}(\overline{\Omega})_\Delta$ when $\mu\in V^{-1,\alpha}(\partial\Omega)$, the validity of the jump formulas for an acoustic single layer with density $\mu\in V^{-1,\alpha}(\partial\Omega)$ and that the integral equation in $V^{-1,\alpha}(\partial\Omega) $ that corresponds to the Neumann problem is in the form of a compact perturbation of the identity in order to apply the Fredholm Alternative Theorem of Wendland \cite{We67}, \cite{We70} in the duality pairing
\[
\left(V^{-1,\alpha}(\partial\Omega),C^{1,\alpha} (\partial\Omega)\right)\,.
\]
We do so by means of Theorems \ref{thm:slhco-1a}, \ref{thm:jufonslap}, and Corollary \ref{corol:thm:wt-1acph}.  

Due to length  restrictions, we plan to develop the consequent  treatment of boundary value problems for the Helmholtz equation with data in $V^{-1,\alpha}(\partial\Omega)$ both in $\Omega$ and in the exterior $\Omega^-$ of $\Omega$     in  forthcoming papers (see \cite{La25b}). 

The paper is organized as follows. Section \ref{sec:prelnot} is a section of
 preliminaries and notation. 
  In section \ref{sec:doun} we introduce the distributional form of the outward unit normal derivative for functions which are locally of class $C^{0,\alpha}(\overline{\Omega^-})_\Delta$ of \cite{La25}. In section \ref{sec:acvolpot} we introduce some   known properties of the acoustic volume potential. In section \ref{sec:acsilah-1a}, we prove Theorem \ref{thm:slhco-1a} on the membership in $C^{0,\alpha}(\overline{\Omega})_\Delta$ of an acoustic single layer potential with density in the space  $V^{-1,\alpha}(\partial\Omega)$ of distributions on the boundary. In section \ref{sec:silahjump}, we prove Theorem \ref{thm:jufonslap} on the validity of the jump formulas of the normal derivative of an acoustic single layer potential with density in  the space  $V^{-1,\alpha}(\partial\Omega)$. In section \ref{sec:qsph} we prove the variant of the Plemelj symmetrization principle of Lemma \ref{lem:qldcomvw}  that involves both harmonic and acoustic  layer potentials. In section \ref{sec:cowht-1a}, we prove Theorem \ref{thm:wt-1ach} and the following Corollary \ref{corol:thm:wt-1acph} that implies that the Fredholm integral equation for  the Neumann problem is in the form of a compact perturbation of the identity.   In the appendix at the end of the paper, we have included a formulation of the classical third  Green	  Identity that we need in the paper.

  \section{Preliminaries and notation}\label{sec:prelnot} Unless otherwise specified,  we assume  throughout the paper that
\[
n\in {\mathbb{N}}\setminus\{0,1\}\,,
\]
where ${\mathbb{N}}$ denotes the set of natural numbers including $0$. 
If $X$ and $Y$, $Z$ are normed spaces, then ${\mathcal{L}}(X,Y)$ denotes the space of linear and continuous maps from $X$ to $Y$ and ${\mathcal{L}}^{(2)}(X\times Y, Z)$ 
denotes the space of bilinear and continuous maps from $X\times Y$ to $Z$ with their usual operator norm (cf.~\textit{e.g.}, \cite[pp.~16, 621]{DaLaMu21}). 
 
$|A|$ denotes the operator norm of a matrix $A$ with real (or complex) entries, 
       $A^{t}$ denotes the transpose matrix of $A$.	 
 
 Let $\Omega$ be an open subset of ${\mathbb{R}}^n$. $C^{1}(\Omega)$ denotes the set of continuously differentiable functions from $\Omega$ to ${\mathbb{C}}$. 
 Let $s\in {\mathbb{N}}\setminus\{0\}$, $f\in \left(C^{1}(\Omega)\right)^{s} $. Then   $Df$ denotes the Jacobian matrix of $f$.   
 
 For the (classical) definition of   open Lipschitz subset of ${\mathbb{R}}^n$ and of   open subset of ${\mathbb{R}}^n$ 
   of class $C^{m}$ or of class $C^{m,\alpha}$
  and of the H\"{o}lder and Schauder spaces $C^{m,\alpha}(\overline{\Omega})$
  on the closure $\overline{\Omega}$ of  an open set $\Omega$ and 
  of the H\"{o}lder and Schauder spaces
   $C^{m,\alpha}(\partial\Omega)$ 
on the boundary $\partial\Omega$ of an open set $\Omega$ for some $m\in{\mathbb{N}}$, $\alpha\in ]0,1]$, we refer for example to
    Dalla Riva, the author and Musolino  \cite[\S 2.3, \S 2.6, \S 2.7, \S 2.9, \S 2.11, \S 2.13,   \S 2.20]{DaLaMu21}.  If $m\in {\mathbb{N}}$, 
 $C^{m}_b(\overline{\Omega})$ denotes the space of $m$-times continuously differentiable functions from $\Omega$ to ${\mathbb{C}}$ such that all 
the partial derivatives up to order $m$ have a bounded continuous extension to    $\overline{\Omega}$ and we set
\[
\|f\|_{   C^{m}_{b}(
\overline{\Omega} )   }\equiv
\sum_{|\eta|\leq m}\, \sup_{x\in \overline{\Omega}}|D^{\eta}f(x)|
\qquad\forall f\in C^{m}_{b}(
\overline{\Omega} )\,,
\]
 where $|\eta|\equiv
\eta_{1}+\dots +\eta_{n}  $ for all $\eta\in {\mathbb{N}}^{n}$.	 
If $\alpha\in ]0,1]$, then 
$C^{m,\alpha}_b(\overline{\Omega})$ denotes the space of functions of $C^{m}_{b}(
\overline{\Omega}) $  such that the  partial derivatives of order $m$ are $\alpha$-H\"{o}lder continuous in $\Omega$. Then we equip $C^{m,\alpha}_{b}(\overline{\Omega})$ with the norm
\[
\|f\|_{  C^{m,\alpha}_{b}(\overline{\Omega})  }\equiv 
\|f\|_{  C^{m }_{b}(\overline{\Omega})  }
+\sum_{|\eta|=m}|D^{\eta}f|_{\alpha}\qquad\forall f\in C^{m,\alpha}_{b}(\overline{\Omega})\,,
\]
where $|D^{\eta}f|_{\alpha}$ denotes the $\alpha$-H\"{o}lder constant of the partial derivative $D^{\eta}f$ of order $\eta$ of $f$ in $\Omega$. If $\Omega$ is bounded, we obviously have $C^{m }_{b}(\overline{\Omega})=C^{m } (\overline{\Omega})$ and $C^{m,\alpha}_{b}(\overline{\Omega})=C^{m,\alpha} (\overline{\Omega})$.  Then $C^{m,\alpha}_{{\mathrm{loc}}}(\overline{\Omega }) $  denotes 
the space  of those functions $f\in C^{m}(\overline{\Omega} ) $ such that $f_{|\overline{\Omega'}} $ belongs to $
C^{m,\alpha}(   \overline{ \Omega' }   )$ for all bounded open subsets $\Omega'$ of ${\mathbb{R}}^n$ such that $\overline{\Omega'}\subseteq\overline{\Omega}$. 
The space of complex valued functions of class $C^m$ with compact support in an open set $\Omega$ of ${\mathbb{R}}^n$ is denoted $C^m_c(\Omega)$ and similarly for $C^\infty_c(\Omega)$. We also set ${\mathcal{D}}(\Omega)\equiv C^\infty_c(\Omega)$. Then the dual ${\mathcal{D}}'(\Omega)$ is known to be the space of distributions in $\Omega$. The support of a function or of a distribution is denoted by the abbreviation `${\mathrm{supp}}$'.  We also set
\[
\Omega^-\equiv {\mathbb{R}}^n\setminus\overline{\Omega}\,.
\] 
We denote by $\nu_\Omega$ or simply by $\nu$ the outward unit normal of $\Omega$ on $\partial\Omega$. Then $\nu_{\Omega^-}=-\nu_\Omega$ is the outward unit normal of $\Omega^-$ on $\partial\Omega=\partial\Omega^-$.
We now summarize the definition and some elementary properties of the Schau\-der space $C^{-1,\alpha}(\overline{\Omega})$ 
by following the presentation of Dalla Riva, the author and Musolino \cite[\S 2.22]{DaLaMu21}.
\begin{definition} 
\label{defn:sch-1}\index{Schauder space!with negative exponent}
 Let $\alpha\in]0,1]$. Let $\Omega$ be a bounded open subset of ${\mathbb{R}}^{n}$. We denote by $C^{-1,\alpha}(\overline{\Omega})$ the subspace 
 \[
 \left\{
 f_{0}+\sum_{j=1}^{n}\frac{\partial}{\partial x_{j}}f_{j}:\,f_{j}\in 
 C^{0,\alpha}(\overline{\Omega})\ \forall j\in\{0,\dots,n\}
 \right\}\,,
 \]
 of the space of distributions ${\mathcal{D}}'(\Omega)$  in $\Omega$ and we set
 \begin{eqnarray}
\label{defn:sch-2}
\lefteqn{
\|f\|_{  C^{-1,\alpha}(\overline{\Omega})  }
\equiv\inf\biggl\{\biggr.
\sum_{j=0}^{n}\|f_{j}\|_{ C^{0,\alpha}(\overline{\Omega})  }
:\,
}
\\ \nonumber
&&\qquad\qquad\qquad\qquad
f=f_{0}+\sum_{j=1}^{n}\frac{\partial}{\partial x_{j}}f_{j}\,,\ 
f_{j}\in C^{0,\alpha}(\overline{\Omega})\ \forall j\in \{0,\dots,n\}
\biggl.\biggr\}\,.
\end{eqnarray}
\end{definition}
$(C^{-1,\alpha}(\overline{\Omega}), \|\cdot\|_{  C^{-1,\alpha}(\overline{\Omega})  })$ is known to be a Banach space and  is continuously embedded into ${\mathcal{D}}'(\Omega)$. Also, the definition of the norm $\|\cdot\|_{  C^{-1,\alpha}(\overline{\Omega})  }$ implies that $C^{0,\alpha}(\overline{\Omega})$ is continuously embedded into $C^{-1,\alpha}(\overline{\Omega})$ and that the partial differentiation $\frac{\partial}{\partial x_{j}}$ is continuous from 
$C^{0,\alpha}(\overline{\Omega})$ to $C^{-1,\alpha}(\overline{\Omega})$ for all $j\in\{1,\dots,n\}$.  Generically, the  elements of $C^{-1,\alpha}(\overline{\Omega})$ are not integrable functions, but distributions in $\Omega$.  Then we have the following statement of \cite[Prop.~3.1]{La24c} that shows that the elements of $C^{-1,\alpha}(\overline{\Omega}) $, which have been defined as elements of	  the dual of ${\mathcal{D}}(\Omega)$, can be extended to elements of the dual of $C^{1,\alpha}(\overline{\Omega})$.  
\begin{proposition}\label{prop:nschext}
 Let $\alpha\in]0,1[$. Let $\Omega$ be a bounded open Lipschitz subset of ${\mathbb{R}}^{n}$.  There exists one and only one  linear and continuous extension operator $E^\sharp_\Omega$ from $C^{-1,\alpha}(\overline{\Omega})$ to $\left(C^{1,\alpha}(\overline{\Omega})\right)'$ such that
 \begin{eqnarray}\label{prop:nschext2}
\lefteqn{
\langle E^\sharp_\Omega[f],v\rangle 
}
\\ \nonumber
&&\ \
 =
\int_{\Omega}f_{0}v\,dx+\int_{\partial\Omega}\sum_{j=1}^{n} (\nu_{\Omega})_{j}f_{j}v\,d\sigma
 -\sum_{j=1}^{n}\int_{\Omega}f_{j}\frac{\partial v}{\partial x_j}\,dx
\quad \forall v\in C^{1,\alpha}(\overline{\Omega})
\end{eqnarray}
for all $f=  f_{0}+\sum_{j=1}^{n}\frac{\partial}{\partial x_{j}}f_{j}\in C^{-1,\alpha}(\overline{\Omega}) $. Moreover, 
\begin{equation}\label{prop:nschext1}
E^\sharp_\Omega[f]_{|\Omega}=f\,, \ i.e.,\ 
\langle E^\sharp_\Omega[f],v\rangle =\langle f,v\rangle \qquad\forall v\in {\mathcal{D}}(\Omega)
\end{equation}
for all $f\in C^{-1,\alpha}(\overline{\Omega})$ and
\begin{equation}\label{prop:nschext3}
\langle E^\sharp_\Omega[f],v\rangle =\langle f,v\rangle \qquad\forall v\in C^{1,\alpha}(\overline{\Omega})
\end{equation}
for all $f\in C^{0,\alpha}(\overline{\Omega})$.
\end{proposition}
When no ambiguity can arise, we simply write $E^\sharp$ instead of $E^\sharp_\Omega$. To see why the extension operator $E^\sharp$ can be considered as `canonical', we refer to \cite[Prop.~7]{La24d}. Next we introduce the following multiplication lemma. For a proof, we refer to \cite[Lem.~2.3]{La25}.
\begin{lemma}\label{lem:multc1ac-1a}
  Let   $\alpha\in ]0,1[$. Let $\Omega$ be  
   a bounded open Lipschitz subset of ${\mathbb{R}}^{n}$.  Then the pointwise product is bilinear and continuous from 
  $C^{1,\alpha}(\overline{\Omega})\times C^{-1,\alpha}(\overline{\Omega})$ to $C^{-1,\alpha}(\overline{\Omega})$. 
\end{lemma}
 Next we introduce our space for the solutions.

\begin{definition}\label{defn:c0ade}
 Let   $\alpha\in ]0,1]$. Let $\Omega$ be a bounded open  subset of ${\mathbb{R}}^{n}$. Let
 \begin{eqnarray}\label{defn:c0ade1}
C^{0,\alpha}(\overline{\Omega})_\Delta
&\equiv&\biggl\{u\in C^{0,\alpha}(\overline{\Omega}):\,\Delta u\in C^{-1,\alpha}(\overline{\Omega})\biggr\}\,,
\\ \nonumber
\|u\|_{ C^{0,\alpha}(\overline{\Omega})_\Delta }
&\equiv& \|u\|_{ C^{0,\alpha}(\overline{\Omega})}
+\|\Delta u\|_{C^{-1,\alpha}(\overline{\Omega})}
\qquad\forall u\in C^{0,\alpha}(\overline{\Omega})_\Delta\,.
\end{eqnarray}
\end{definition}
Since $C^{0,\alpha}(\overline{\Omega})$ and $C^{-1,\alpha}(\overline{\Omega})$ are Banach spaces,   $\left(\|u\|_{ C^{0,\alpha}(\overline{\Omega})_\Delta }, \|\cdot \|_{ C^{0,\alpha}(\overline{\Omega})_\Delta }\right)$ is a Banach space.  For subsets $\Omega$ that are not necessarily bounded, we introduce the following statement of 
\cite[Lem.~2.5]{La25}.
\begin{lemma}\label{lem:c1alcof}
 Let   $\alpha\in ]0,1]$. Let $\Omega$ be an open  subset of ${\mathbb{R}}^{n}$. Then the space
\begin{eqnarray}\label{lem:c1alcof1}
\lefteqn{
 C^{0,\alpha}_{	{\mathrm{loc}}	}(\overline{\Omega})_\Delta\equiv\biggl\{
 f\in C^{0}(\overline{\Omega}):\, f_{|\overline{\Omega'}}\in C^{0,\alpha}(\overline{\Omega'})_\Delta\ \text{for\ all\  } 
 }
\\ \nonumber
&&\qquad\qquad\qquad\qquad 
 \text{bounded\ open\ subsets}\ \Omega'\ \text{of}\ {\mathbb{R}}^n\ \text{such\ that} \  \overline{\Omega'}\subseteq\overline{\Omega}  \biggr\}\,,
\end{eqnarray}
 with the family of seminorms
\begin{eqnarray}\label{lem:c1alcof2}
\lefteqn{
{\mathcal{P}}_{C^{0,\alpha}_{	{\mathrm{loc}}	}(\overline{\Omega})_\Delta}\equiv
\biggl\{
\|\cdot\|_{C^{0,\alpha}(\overline{\Omega'})_\Delta}:\,  
}
\\ \nonumber
&& \qquad\qquad\qquad\qquad 
 \Omega'\ \text{is\ a\ bounded\ open  subset\ of\ }  {\mathbb{R}}^n \ \text{such\ that} \  \overline{\Omega'}\subseteq\overline{\Omega}
\biggr\}
\end{eqnarray}
is a Fr\'{e}chet space.
\end{lemma}

Next we introduce the following multiplication lemma. For a proof, we refer to \cite[Lem.~2.6]{La25}
\begin{lemma}\label{lem:multc2ac0ad}
  Let   $\alpha\in ]0,1[$. Let $\Omega$ be a bounded open Lipschitz subset of ${\mathbb{R}}^{n}$. Then the pointwise product is bilinear and continuous from 
  $C^{2,\alpha}(\overline{\Omega})\times C^{0,\alpha}(\overline{\Omega})_\Delta$ to $C^{0,\alpha}(\overline{\Omega})_\Delta$. 
\end{lemma}
 We also note that the following elementary lemma holds. For a proof, we refer to 
 \cite[Lem.~2.7]{La25}
\begin{lemma}\label{lem:c1alf}
 Let   $\alpha\in ]0,1[$. Let $\Omega$ be an open  subset of ${\mathbb{R}}^{n}$. Then the space
\begin{eqnarray}\label{lem:c1alf1}
\lefteqn{
 C^{1,\alpha}_{	{\mathrm{loc}}	}(\Omega)\equiv\biggl\{
 f\in C^{1}(\Omega):\, f_{|\overline{\Omega'}}\in C^{1,\alpha}(\overline{\Omega'})\ \text{for\ all}\  \Omega'\ \text{such\ that}\ 
 }
\\ \nonumber
&&\qquad\qquad 
 \Omega'\ \text{is\ a\ bounded\ open  subset \ of}\ {\mathbb{R}}^n \,, \overline{\Omega'}\subseteq\Omega  \biggr\}\,,
\end{eqnarray}
 with the family of seminorms
\begin{eqnarray}\label{lem:c1alf2}
\lefteqn{
{\mathcal{P}}_{C^{1,\alpha}_{{\mathrm{loc}}}(\Omega)}\equiv
\biggl\{
\|\cdot\|_{C^{1,\alpha}(\overline{\Omega'})}:\,  
}
\\ \nonumber
&& \qquad\qquad 
 \Omega'\ \text{is\ a\ bounded\ open  \ subset \ of} \ {\mathbb{R}}^n \,, \overline{\Omega'}\subseteq\Omega  
\biggr\}
\end{eqnarray}
is a Fr\'{e}chet space.
\end{lemma}
 Then we can prove the following. For a proof, we refer to \cite[Lem.~2.8]{La25}.
\begin{lemma}\label{lem:c0ademc1a}
  Let   $\alpha\in ]0,1[$. Let $\Omega$ be a bounded open  subset of ${\mathbb{R}}^{n}$. Then $C^{0,\alpha}(\overline{\Omega})_\Delta$ is continuously embedded into the 
  Fr\'{e}chet space $C^{1,\alpha}_{{\mathrm{loc}}}(\Omega)$ with the family of seminorms ${\mathcal{P}}_{C^{1,\alpha}_{{\mathrm{loc}}}(\Omega)}$.
\end{lemma}
 Next we   introduce a subspace of  $C^{1,\alpha}(\overline{\Omega})$ that we need in the sequel.  
\begin{definition}\label{defn:c1ade}
 Let   $\alpha\in ]0,1]$. Let $\Omega$ be a bounded open  subset of ${\mathbb{R}}^{n}$ of class $C^{1,\alpha}$. Let
 \begin{eqnarray}\label{defn:c1ade1}
C^{1,\alpha}(\overline{\Omega})_\Delta
&\equiv&\biggl\{u\in C^{1,\alpha}(\overline{\Omega}):\,\Delta u\in C^{0,\alpha}(\overline{\Omega})\biggr\}\,,
\\ \nonumber
\|u\|_{ C^{1,\alpha}(\overline{\Omega})_\Delta }
&\equiv& \|u\|_{ C^{1,\alpha}(\overline{\Omega})}
+\|\Delta u\|_{C^{0,\alpha}(\overline{\Omega})}
\qquad\forall u\in C^{1,\alpha}(\overline{\Omega})_\Delta\,.
\end{eqnarray}
\end{definition}
If $\Omega$ is a bounded open subset of ${\mathbb{R}}^{n}$, then 
 $C^{1,\alpha}(\overline{\Omega})$ and $C^{0,\alpha}(\overline{\Omega})$ are Banach spaces and accordingly   $\left(\|u\|_{ C^{1,\alpha}(\overline{\Omega})_\Delta }, \|\cdot \|_{ C^{1,\alpha}(\overline{\Omega})_\Delta }\right)$ is a Banach space.  Moreover, by applying \cite[Lem.~5.12]{La24c} to each closed ball  contained in $\Omega$, we deduce that
\begin{equation}\label{eq:c1adc2}
C^{1,\alpha}(\overline{\Omega})_\Delta\subseteq C^2(\Omega)\,.
\end{equation}
 We also note that  if $\Omega$ is a bounded open Lipschitz subset of ${\mathbb{R}}^n$, then 
$C^{1,\alpha}(\overline{\Omega})_\Delta$ is continuously embedded into $ C^{0,\alpha}(\overline{\Omega})_\Delta$.  Next we introduce the following classical result on the Green operator for the interior Dirichlet problem. For a proof, we refer for example to \cite[Thm.~4.8]{La24b}.
 \begin{theorem}\label{thm:idwp}
Let $m\in {\mathbb{N}}$, $\alpha\in ]0,1[$. Let $\Omega$ be a bounded open  subset of ${\mathbb{R}}^{n}$ of class $C^{\max\{m,1\},\alpha}$.
Then the map ${\mathcal{G}}_{\Omega,d,+}$  from $C^{m,\alpha}(\partial\Omega)$ to the closed subspace 
 \begin{equation}\label{thm:idwp1}
C^{m,\alpha}_h(\overline{\Omega}) \equiv \{
u\in C^{m,\alpha}(\overline{\Omega}), u\ \text{is\ harmonic\ in}\ \Omega\}
\end{equation}
of $ C^{m,\alpha}(\overline{\Omega})$ that takes $v$ to the only solution $v^\sharp$ of the Dirichlet problem
\begin{equation}\label{defn:cinspo3}
\left\{
\begin{array}{ll}
 \Delta v^\sharp=0 & \text{in}\ \Omega\,,
 \\
v^\sharp_{|\partial\Omega} =v& \text{on}\ \partial\Omega 
\end{array}
\right.
\end{equation}
is a linear homeomorphism.
\end{theorem}
Next we introduce the following  approximation lemma.
\begin{lemma}\label{lem:apr1a}
 Let $\alpha\in]0,1]$. 
 Let $\Omega$ be a bounded open subset of ${\mathbb{R}}^n$ of class $C^{1,\alpha}$.  If $g\in C^{1,\alpha}(\overline{\Omega})$, then there exists a sequence $\{g_j\}_{j\in {\mathbb{N}}}$ in 
 $C^{\infty}(\overline{\Omega})$ such that
 \begin{equation}\label{lem:apr1a1}
 \sup_{j\in {\mathbb{N}}}\|g_j\|_{
 C^{1,\alpha}(\overline{\Omega}) 
 }<+\infty\,,\qquad
 \lim_{j\to\infty}g_j=g\quad\text{in}\ C^{1,\beta}(\overline{\Omega})\quad \forall\beta\in]0,\alpha[\,.
 \end{equation}
\end{lemma}
For a proof, we refer to \cite[Lem.~A.3]{La24c}.   Next we plan to introduce the normal derivative of the functions in  $C^{0,\alpha}(\overline{\Omega})_\Delta$ as in \cite{La24c}. To do so,    we introduce the (classical) interior Steklov-Poincar\'{e} operator (or interior  Dirichlet-to-Neumann map).
\begin{definition}\label{defn:cinspo}
 Let $\alpha\in]0,1[$.  Let  $\Omega$ be a  bounded open subset of ${\mathbb{R}}^{n}$ of class $C^{1,\alpha}$. The classical interior Steklov-Poincar\'{e} operator is defined to be the operator $S_{\Omega,+}$ from
 \begin{equation}\label{defn:cinspo1}
C^{1,\alpha}(\partial\Omega)\quad\text{to}\quad C^{0,\alpha}(\partial\Omega)
\end{equation}
that takes $v\in C^{1,\alpha}(\partial\Omega)$ to the function 
 \begin{equation}\label{defn:cinspo2}
S_{\Omega,+}[v](x)\equiv \frac{\partial  }{\partial\nu}{\mathcal{G}}_{\Omega,d,+}[v](x)\qquad\forall x\in\partial\Omega\,.
\end{equation}
 \end{definition}
   Since   the classical normal derivative is continuous from $C^{1,\alpha}(\overline{\Omega})$ to $C^{0,\alpha}(\partial\Omega)$, the continuity of ${\mathcal{G}}_{\Omega,d,+}$ implies  that $S_{\Omega,+}[\cdot]$ is linear and continuous from 
  $C^{1,\alpha}(\partial\Omega)$ to $C^{0,\alpha}(\partial\Omega)$. Then we have the following definition of \cite[(41)]{La24c}.
  \begin{definition}\label{defn:conoderdedu}
 Let $\alpha\in]0,1[$.  Let  $\Omega$ be a  bounded open subset of ${\mathbb{R}}^{n}$ of class $C^{1,\alpha}$. If  $u\in C^{0}(\overline{\Omega})$ and $\Delta u\in  C^{-1,\alpha}(\overline{\Omega})$, then we define the distributional  normal derivative $\partial_\nu u$
 of $u$ to be the only element of the dual $(C^{1,\alpha}(\partial\Omega))'$ that satisfies the following equality
 \begin{equation}\label{defn:conoderdedu1}
\langle \partial_\nu u ,v\rangle \equiv\int_{\partial\Omega}uS_{\Omega,+}[v]\,d\sigma
+\langle E^\sharp_\Omega[\Delta u],{\mathcal{G}}_{\Omega,d,+}[v]\rangle 
\qquad\forall v\in C^{1,\alpha}(\partial\Omega)\,.
\end{equation}
\end{definition}
 The normal derivative of Definition \ref{defn:conoderdedu} extends the classical one in the sense that if $u\in C^{1,\alpha}(\overline{\Omega})$, then under the assumptions on $\alpha$ and $\Omega$  of Definition \ref{defn:conoderdedu}, we have 
 \begin{equation}\label{lem:conoderdeducl1}
 \langle \partial_\nu u ,v\rangle \equiv\int_{\partial\Omega}\frac{\partial u}{\partial\nu}v\,d\sigma
  \quad\forall v\in C^{1,\alpha}(\partial\Omega)\,,
\end{equation}
where $\frac{\partial u}{\partial\nu}$ in the right hand side denotes the classical normal derivative of $u$ on $\partial\Omega$ (cf.~\cite[Lem.~5.5]{La24c}). In the sequel, we use the classical symbol $\frac{\partial u}{\partial\nu}$ also for  $\partial_\nu u$ when no ambiguity can arise.\par  
  
Next we introduce the  function space $V^{-1,\alpha}(\partial\Omega)$ on the boundary of $\Omega$ for the normal derivatives of the functions of $C^{0,\alpha}(\overline{\Omega})_\Delta$ as in  \cite[Defn.~13.2, 15.10, Thm.~18.1]{La24b}.
\begin{definition}\label{defn:v-1a}
Let   $\alpha\in ]0,1[$. Let $\Omega$ be a bounded open  subset of ${\mathbb{R}}^{n}$ of class $C^{1,\alpha}$. Let 
\begin{eqnarray}\label{defn:v-1a1}
 \lefteqn{V^{-1,\alpha}(\partial\Omega)\equiv \biggl\{\mu_0+S_{\Omega,+}^t[\mu_1]:\,\mu_0, \mu_1\in C^{0,\alpha}(\partial\Omega)
\biggr\}\,,
}
\\ \nonumber
\lefteqn{
\|\tau\|_{  V^{-1,\alpha}(\partial\Omega) }
\equiv\inf\biggl\{\biggr.
 \|\mu_0\|_{ C^{0,\alpha}(\partial\Omega)  }+\|\mu_1\|_{ C^{0,\alpha}(\partial\Omega)  }
:\,
 \tau=\mu_0+S_{\Omega,+}^t[\mu_1]\biggl.\biggr\}\,,
 }
 \\ \nonumber
 &&\qquad\qquad\qquad\qquad\qquad\qquad\qquad\qquad\qquad
 \forall \tau\in  V^{-1,\alpha}(\partial\Omega)\,,
\end{eqnarray}
where $S_{\Omega,+}^t$ is the transpose map of $S_{\Omega,+}$.
\end{definition}
As shown in \cite[\S 13]{La24b},  $(V^{-1,\alpha}(\partial\Omega), \|\cdot\|_{  V^{-1,\alpha}(\partial\Omega)  })$ is a Banach space that is continuously embedded into the dual of $C^{1,\alpha}(\partial\Omega)$. By definition of the norm, $C^{0,\alpha}(\partial\Omega)$ is continuously embedded into $V^{-1,\alpha}(\partial\Omega)$. Moreover, we have the following statement  on the continuity of the normal derivative on $C^{0,\alpha}(\overline{\Omega})_\Delta$  (see  \cite[Prop.~6.6]{La24c}). 
\begin{proposition}\label{prop:ricodnu}
 Let   $\alpha\in ]0,1[$. Let $\Omega$ be a bounded open  subset of 
 ${\mathbb{R}}^{n}$ of class $C^{1,\alpha}$. Then the distributional normal derivative   $\partial_\nu$ is a continuous surjection of $C^{0,\alpha}(\overline{\Omega})_\Delta$ onto $V^{-1,\alpha}(\partial\Omega)$ and there exists $Z\in {\mathcal{L}}\left(V^{-1,\alpha}(\partial\Omega),C^{0,\alpha}(\overline{\Omega})_\Delta\right)$ such that
 \begin{equation}\label{prop:ricodnu1}
\partial_\nu Z[g]=g\qquad\forall g\in V^{-1,\alpha}(\partial\Omega)\,,
\end{equation}
\textit{i.e.}, $Z$ is a right inverse of  $\partial_\nu$.
\end{proposition}
We also mention that the  Definition \ref{defn:conoderdedu} of normal derivative admits the following	 	 different formulation at least in case $u\in C^{0,\alpha}(\overline{\Omega})_\Delta$ (see \cite[Prop.~5.15]{La24c}).
\begin{proposition}\label{prop:node1adeq}
 Let   $\alpha\in ]0,1[$. Let $\Omega$ be a bounded open  subset of 
 ${\mathbb{R}}^{n}$ of class $C^{1,\alpha}$. Let $E_\Omega$ be a linear map from $C^{1,\alpha}(\partial\Omega)$ to $ C^{1,\alpha}(\overline{\Omega})_\Delta$ such that
 \begin{equation}\label{prop:node1adeq0}
 E_\Omega[f]_{|\partial\Omega}=f\qquad\forall f\in C^{1,\alpha}(\partial\Omega)\,.
 \end{equation}
 If $u\in C^{0,\alpha}(\overline{\Omega})_\Delta$, then the distributional  normal derivative $\partial_\nu u$
 of $u$  is characterized by the validity of the following equality
 \begin{eqnarray}\label{prop:node1adeq1}
 \lefteqn{
\langle \partial_\nu u ,v\rangle =\int_{\partial\Omega}u
\frac{\partial}{\partial\nu}E_\Omega[v]
\,d\sigma
}
\\ \nonumber
&&\qquad\qquad
+\langle E^\sharp_\Omega[\Delta u],E_\Omega[v]\rangle -\int_\Omega\Delta (E_\Omega[v]) u\,dx
\qquad\forall v\in C^{1,\alpha}(\partial\Omega)\,.
\end{eqnarray}
\end{proposition}
   We note that in equality (\ref{prop:node1adeq1}) we can use the `extension' operator $E_\Omega$ that we prefer and that accordingly equality (\ref{prop:node1adeq1}) is independent of the specific choice of $E_\Omega$. When we deal with problems for the Laplace operator, a good choice is
 $E_\Omega={\mathcal{G}}_{\Omega,d,+}$, so that the last term in the right hand side of (\ref{prop:node1adeq1}) disappears.  
 In particular, under the assumptions of Proposition \ref{prop:node1adeq}, an extension operator as $E_\Omega$ always exists.
 
 \section{A  distributional outward unit normal derivative for the functions of $C^{0,\alpha}_{
{\mathrm{loc}}	}(\overline{\Omega^-})_\Delta$}
\label{sec:doun}
 
 We now plan to define the normal derivative on the boundary the functions of $C^{0,\alpha}_{
{\mathrm{loc}}	}(\overline{\Omega^-})_\Delta$ with $\alpha\in]0,1[$ in case $\Omega$ is a bounded open subset of ${\mathbb{R}}^n$ of class $C^{1,\alpha}$ as in \cite[\S 3]{La25}.  To do so, we  choose $r\in]0,+\infty[$ such that $\overline{\Omega}\subseteq {\mathbb{B}}_n(0,r)$ and we observe that the linear operator from $C^{1,\alpha}(\partial\Omega)$ to $C^{1,\alpha}((\partial\Omega)\cup(\partial{\mathbb{B}}_n(0,r)))$ that is defined by the equality
\[
\stackrel{o}{E}_{\partial\Omega,r}[v]\equiv\left\{
\begin{array}{ll}
 v(x) & \text{if}\ x\in \partial\Omega\,,
 \\
 0 & \text{if}\ x\in \partial{\mathbb{B}}_n(0,r)\,,
\end{array}
\right.\qquad\forall v\in C^{1,\alpha}(\partial\Omega) 
\]
is linear and continuous and that accordingly the transpose map $\stackrel{o}{E}_{\partial\Omega,r}^t$ is linear and continuous
\[
\text{from}\ 
\left(C^{1,\alpha}((\partial\Omega)\cup(\partial{\mathbb{B}}_n(0,r)))\right)'\ 
\text{to}\ \left(C^{1,\alpha}(\partial\Omega)\right)'\,.
\]
Then we introduce the following definition as in \cite[\S 3]{La25}.
\begin{definition}\label{defn:endedr}
 Let   $\alpha\in ]0,1[$. Let $\Omega$ be a bounded open  subset of 
 ${\mathbb{R}}^{n}$ of class $C^{1,\alpha}$. Let $r\in]0,+\infty[$ such that $\overline{\Omega}\subseteq {\mathbb{B}}_n(0,r)$. If $u\in C^{0,\alpha}_{
{\mathrm{loc}}	}(\overline{\Omega^-})_\Delta$, then we set
\begin{equation}\label{defn:endedr1}
 \partial_{\nu_{\Omega^-}}u
 \equiv 
 \stackrel{o}{E}_{\partial\Omega,r}^t\left[
\partial_{\nu_{{\mathbb{B}}_n(0,r)\setminus\overline{\Omega}}}u
\right] \,.	 
\end{equation}
 \end{definition}
 Under the assumptions of Definition \ref{defn:endedr},   $\partial_{\nu_{\Omega^-}}u$ is an element of $\left(C^{1,\alpha}(\partial\Omega)\right)'$. Under the assumptions of Definition \ref{defn:endedr}, there exists  a linear (extension) map $E_{{\mathbb{B}}_n(0,r)\setminus\overline{\Omega}}$ from $C^{1,\alpha}((\partial\Omega)\cup(\partial
{\mathbb{B}}_n(0,r)
))$ to $ C^{1,\alpha}(\overline{{\mathbb{B}}_n(0,r)}\setminus\Omega)_\Delta$ such that
 \begin{equation}\label{prop:node1adeq0r}
 E_{{\mathbb{B}}_n(0,r)\setminus\overline{\Omega}}[f]_{|(\partial\Omega)\cup(\partial
{\mathbb{B}}_n(0,r)
)}=f\qquad\forall f\in C^{1,\alpha}((\partial\Omega)\cup(\partial
{\mathbb{B}}_n(0,r)
))\,.
 \end{equation}
Then the definition of normal derivative on the boundary of 
${\mathbb{B}}_n(0,r)\setminus\overline{\Omega}$ implies that
\begin{eqnarray}\label{defn:endedr2}
\lefteqn{
\langle \partial_{\nu_{\Omega^-}}u,v\rangle  
=-\int_{\partial\Omega}u\frac{\partial}{\partial\nu_{\Omega}}E_{{\mathbb{B}}_n(0,r)\setminus\overline{\Omega}}[\stackrel{o}{E}_{\partial\Omega,r}[v]]\,d\sigma
}
\\ \nonumber
&&\qquad\qquad 
+\int_{\partial{\mathbb{B}}_n(0,r)}u\frac{\partial}{\partial\nu_{{\mathbb{B}}_n(0,r)}}E_{{\mathbb{B}}_n(0,r)\setminus\overline{\Omega}}[\stackrel{o}{E}_{\partial\Omega,r}[v]]\,d\sigma
\\ \nonumber
&&\qquad\qquad 
+\langle E^\sharp_{{\mathbb{B}}_n(0,r)\setminus\overline{\Omega}}[\Delta u],E_{{\mathbb{B}}_n(0,r)\setminus\overline{\Omega}}[\stackrel{o}{E}_{\partial\Omega,r}[v]]\rangle 
\\ \nonumber
&&\qquad\qquad 
-\int_{{\mathbb{B}}_n(0,r)\setminus\overline{\Omega}}\Delta (E_{{\mathbb{B}}_n(0,r)\setminus\overline{\Omega}}[\stackrel{o}{E}_{\partial\Omega,r}[v]]) u\,dx
\qquad\forall v\in C^{1,\alpha}(\partial\Omega)\,,
\end{eqnarray}
(cf.~(\ref{prop:node1adeq1})).   As shown in \cite[\S 3]{La25},  Definition  \ref{defn:endedr} is independent of the choice of $r\in]0,+\infty[$ such that $\overline{\Omega}\subseteq {\mathbb{B}}_n(0,r)$. Then we also state the following remark of  \cite[\S 3]{La25}.
\begin{remark}\label{rem:conoderdeducl-}
 Let   $\alpha\in ]0,1[$. Let $\Omega$ be a bounded open  subset of 
 ${\mathbb{R}}^{n}$ of class $C^{1,\alpha}$. Let $u\in C^{1,\alpha}_{
{\mathrm{loc}}	}(\overline{\Omega^-})$. Then 
 \begin{equation}\label{rem:conoderdeducl-1}
 \langle \partial_{\nu_{\Omega^-}} u ,v\rangle \equiv\int_{\partial\Omega}\frac{\partial u}{\partial\nu_{\Omega^-}}v\,d\sigma
  \quad\forall v\in C^{1,\alpha}(\partial\Omega)\,,
\end{equation}
where $\frac{\partial u}{\partial\nu_{\Omega^-}}$ in the right hand side denotes the classical $\nu_{\Omega^-}$-normal derivative of $u$ on $\partial\Omega$.  In the sequel, we use the classical symbol $\frac{\partial u}{\partial\nu_{\Omega^-}}$ also for  $\partial_{\nu_{\Omega^-}} u$ when no ambiguity can arise.\par 
\end{remark}

Next we note that if  $u\in C^{0,\alpha}_{{\mathrm{loc}}}(\overline{\Omega^-})$ and $u$ is both harmonic in $\Omega^-$ and harmonic at infinity, then $u\in C^{0,\alpha}_{
{\mathrm{loc}}	}(\overline{\Omega^-})_\Delta$. As shown in \cite[\S 3]{La25}, 
  the distribution $\partial_{\nu_{\Omega^-}}u$ of Definition \ref{defn:endedr} coincides with the normal derivative that has been introduced in \cite[Defn.~6.4]{La24b} for harmonic functions.   By \cite[Thm.~1.1]{La25}, we have the following continuity statement for the distributional normal derivative $\partial_{\nu_{\Omega^-}}$ of
Definition \ref{defn:endedr}.

\begin{proposition}\label{prop:recodnu}
 Let   $\alpha\in ]0,1[$. Let $\Omega$ be a bounded open  subset of 
 ${\mathbb{R}}^{n}$ of class $C^{1,\alpha}$.  Then the distributional normal derivative   $\partial_{\nu_{\Omega^-}}$ is a continuous surjection of $C^{0,\alpha}_{
{\mathrm{loc}}	}(\overline{\Omega^-})_\Delta$ onto $V^{-1,\alpha}(\partial\Omega)$ and there exists a linear and continuous  operator   $Z_-$ from  $ V^{-1,\alpha}(\partial\Omega)$ to $C^{0,\alpha}_{
{\mathrm{loc}}	}(\overline{\Omega^-})_\Delta $ such that
 \begin{equation}\label{prop:recodnu1}
\partial_{\nu_{\Omega^-}} Z_-[g]=g\qquad\forall g\in V^{-1,\alpha}(\partial\Omega)\,,
\end{equation}
\textit{i.e.}, $Z_-$ is a right inverse of  $\partial_{\nu_{\Omega^-}}$. (See Lemma \ref{lem:c1alcof} for the topology of $C^{0,\alpha}_{
{\mathrm{loc}}	}(\overline{\Omega^-})_\Delta $).
\end{proposition}
 
\section{Preliminaries to the acoustic volume potential}
 \label{sec:acvolpot}
Let $\alpha\in]0,1]$ and $m\in {\mathbb{N}}$.  If $\Omega$ is a  bounded open subset of ${\mathbb{R}}^{n}$, then we can consider the restriction map $r_{|\overline{\Omega}}$ from ${\mathcal{D}}({\mathbb{R}}^n)$ to $C^{m,\alpha}(\overline{\Omega})$. Then the transpose map $r_{|\overline{\Omega}}^t$ is linear and continuous from $(C^{m,\alpha}(\overline{\Omega}))'$ to ${\mathcal{D}}'({\mathbb{R}}^n)$. Moreover, if $\mu\in (C^{m,\alpha}(\overline{\Omega}))'$, then $r_{|\overline{\Omega}}^t\mu$ has compact support. Hence, it makes sense to consider the convolution of 
  $r_{|\overline{\Omega}}^t\mu$ with  the fundamental solution of either the Laplace or the Helmholtz  operator. Thus we are now ready to introduce the following known definition.
 \begin{definition}\label{defn:dvpsl}
 Let $\alpha\in]0,1]$, $m\in {\mathbb{N}}$. 
 Let   $\Omega$ be a bounded open subset of ${\mathbb{R}}^{n}$. Let $\lambda\in {\mathbb{C}}$.  Let $S_{n,\lambda} $ be a fundamental solution 
of the operator $\Delta+\lambda$. If $\mu\in (C^{m,\alpha}(\overline{\Omega}))'$, then the (distributional) volume potential relative to $S_{n,\lambda} $ and $\mu$ is the distribution
\[
{\mathcal{P}}_\Omega[S_{n,\lambda} ,\mu]=(r_{|\overline{\Omega}}^t\mu)\ast S_{n,\lambda}  \in {\mathcal{D}}'({\mathbb{R}}^n)\,.
\]
\end{definition}
Under the assumptions of Definition \ref{defn:dvpsl}, we  set
\begin{eqnarray}\label{prop:dvpsl3}
{\mathcal{P}}_\Omega^+[S_{n,\lambda} ,\mu]&\equiv&\left((r_{|\overline{\Omega}}^t\mu)\ast S_{n,\lambda}  \right)_{|\Omega}
\qquad\text{in}\ \Omega\,,
\\ \nonumber
{\mathcal{P}}_\Omega^-[S_{n,\lambda} ,\mu]  &\equiv&
\left((r_{|\overline{\Omega}}^t\mu)\ast S_{n,\lambda}  \right)_{|\Omega^-}
\qquad\text{in}\ \Omega^-\,.
\end{eqnarray}
In general, $(r_{|\overline{\Omega}}^t\mu)\ast S_{n,\lambda} $ is not a function, \textit{i.e.} $(r_{|\overline{\Omega}}^t\mu)\ast S_{n,\lambda} $ is not a distribution that is associated to a locally integrable function in ${\mathbb{R}}^n$. 
However, this is the case if for example  $\mu$ is associated to a function $f$ of $ L^\infty(\Omega)$,  and thus 
 the (distributional) volume potential relative to $S_{n,\lambda} $ and $\mu$ is associated to the function
\begin{equation}\label{prop:dvpsa1}
\int_{\Omega}S_{n,\lambda} (x-y)f(y)\,dy\qquad{\mathrm{a.a.}}\ x\in {\mathbb{R}}^n\,,
\end{equation}
that is locally integrable in ${\mathbb{R}}^n$  and that with some abuse of notation we  denote by the symbol    ${\mathcal{P}}_\Omega[S_{n,\lambda} , f]$.  Then the following classical result is known (cf.~\textit{e.g.}, 
\cite[Thm.~20]{La24d}).
\begin{theorem}\label{thm:nwtdma} 
 Let $m\in {\mathbb{N}}$, $\alpha\in]0,1[$. Let $\Omega$ be a bounded open subset of  ${\mathbb{R}}^n$ of class $C^{m+1,\alpha}$. Let $\lambda\in {\mathbb{C}}$.  Let $S_{n,\lambda}  $ be a fundamental solution 
of the operator $\Delta+\lambda$. Then the following statements hold. 
 
 \item[(i)]  ${\mathcal{P}}_\Omega^+[S_{n,\lambda} ,\cdot]$ is linear and continuous from $C^{m,\alpha}(\overline{\Omega})$ to $C^{m+2,\alpha}(\overline{\Omega})$.
\item[(ii)]   ${\mathcal{P}}_\Omega^-[S_{n,\lambda} ,\cdot]_{|\overline{{\mathbb{B}}_n(0,r)}\setminus\Omega}$ is linear and continuous from the space  $C^{m,\alpha}(\overline{\Omega})$ to the space    $C^{m+2,\alpha}(\overline{{\mathbb{B}}_n(0,r)}\setminus\Omega)$ for all $r\in]0,+\infty[$ such that $\overline{\Omega}\subseteq {\mathbb{B}}_n(0,r)$.
\end{theorem}
 Instead, for a Schauder space of negative exponent, the following statement holds (cf.~\textit{e.g.}, \cite[Prop.~23]{La24d}), that is a  generalization to volume potentials of  nonhomogeneous second order elliptic operators of a known result for the Laplace operator  (cf.~\cite[Thm.~3.6 (ii)]{La08a}, Dalla Riva, the author and Musolino~\cite[Thm.~7.19]{DaLaMu21}).
    \begin{proposition}\label{prop:dvpsnecr-1a}
  Let $\alpha\in]0,1[$.    Let $\Omega$ be a bounded open subset of ${\mathbb{R}}^{n}$ of class $C^{1,\alpha}$. Let $\lambda\in {\mathbb{C}}$.  Let $S_{n,\lambda}  $ be a fundamental solution 
of the operator $\Delta+\lambda$.   Then the following statements hold. 
 \begin{enumerate}
\item[(i)] If  $f=  f_{0}+\sum_{j=1}^{n}\frac{\partial}{\partial x_{j}}f_{j}\in C^{-1,\alpha}(\overline{\Omega}) $, then
 \begin{equation}\label{prop:dvpsnecr-1a2}
{\mathcal{P}}_\Omega^+[S_{n,\lambda} ,E^\sharp[f]]\in C^{1,\alpha}(\overline{\Omega}), \ 
{\mathcal{P}}_\Omega^-[S_{n,\lambda} ,E^\sharp[f]]\in C^{1,\alpha}_{{\mathrm{loc}} }(\overline{\Omega^-}) \,,\end{equation}
and 
\begin{equation}\label{defn:Ppm1}
{\mathcal{P}}_\Omega^+[S_{n,\lambda} ,E^\sharp[f]](x)={\mathcal{P}}_\Omega^-[S_{n,\lambda} ,E^\sharp[f]](x)\qquad\forall x\in\partial\Omega\,.
\end{equation}
Moreover,
\begin{eqnarray}\label{prop:dvpsnecr-1a2a}
&&\Delta {\mathcal{P}}_\Omega^+[S_{n,\lambda} ,E^\sharp[f]] 
+\lambda{\mathcal{P}}_\Omega^+[S_{n,\lambda} ,E^\sharp[f]]= f\qquad\textit{in}\ {\mathcal{D}}'(\Omega)\,,
\\ \nonumber
&&\Delta  {\mathcal{P}}_\Omega^-[S_{n,\lambda} ,E^\sharp[f]]
+\lambda  {\mathcal{P}}_\Omega^-[S_{n,\lambda} ,E^\sharp[f]]= 0\qquad\textit{in}\ {\mathcal{D}}'({\mathbb{R}}^n\setminus\overline{\Omega})
\,.
\end{eqnarray}
\item[(ii)] The linear operator  ${\mathcal{P}}_\Omega^+[S_{n,\lambda} ,E^\sharp[\cdot]]$ is  continuous from $C^{-1,\alpha}(\overline{\Omega}) $ to $C^{1,\alpha}(\overline{\Omega})$.
\item[(iii)] Let $r\in ]0,+\infty[$ be such that $\overline{\Omega}\subseteq {\mathbb{B}}_n(0,r)$. Then  ${\mathcal{P}}_\Omega^-[S_{n,\lambda} ,E^\sharp[\cdot]]_{|\overline{{\mathbb{B}}_n(0,r)}\setminus\Omega}$ is linear and continuous from $C^{-1,\alpha}(\overline{\Omega}) $ to  $C^{1,\alpha}(\overline{{\mathbb{B}}_n(0,r)}\setminus\Omega)$.
\end{enumerate}  
\end{proposition}

\section{Acoustic layer potentials with densities in the space $V^{-1,\alpha}(\partial\Omega)$}
 \label{sec:acsilah-1a}
 
 Let  $\alpha\in]0,1]$. Let  $\Omega$ be a bounded open subset of ${\mathbb{R}}^{n}$ of class $C^{1,\alpha}$. Let  $r_{|\partial\Omega}$  be the restriction map  from ${\mathcal{D}}({\mathbb{R}}^n)$ to $C^{1,\alpha}(\partial\Omega)$. Let $\lambda\in {\mathbb{C}}$. If $S_{n,\lambda} $ is a fundamental solution 
of the operator $\Delta+\lambda$ and $\mu\in (C^{1,\alpha}(\partial\Omega))'$, then the  (distributional) single layer  potential relative to $S_{n,\lambda} $ and $\mu$ is the distribution
\[
v_\Omega[S_{n,\lambda} ,\mu]=(r_{|\partial\Omega}^t\mu)\ast S_{n,\lambda}  \in {\mathcal{D}}'({\mathbb{R}}^n) 
\]
 and we  set
\begin{eqnarray}\label{eq:dsila}
v_\Omega^+[S_{n,\lambda} ,\mu]&\equiv&\left((r_{|\partial\Omega}^t\mu)\ast S_{n,\lambda}  \right)_{|\Omega}
\qquad\text{in}\ \Omega\,,
\\ \nonumber
v_\Omega^-[S_{n,\lambda} ,\mu] &\equiv&
\left((r_{|\partial\Omega}^t\mu)\ast S_{n,\lambda}  \right)_{|\Omega^-}
\qquad\text{in}\ \Omega^-\,.
\end{eqnarray}
It is also known that the restriction of $v_\Omega[S_{n,\lambda} ,\mu]$ to ${\mathbb{R}}^n\setminus\partial\Omega$ equals the (distribution that is associated to) the function
\[
 \langle (r_{|\partial\Omega}^t\mu)(y),S_{n,\lambda} (\cdot-y)\rangle \,.
\]
In the cases in which both $v_\Omega^+[S_{n,\lambda} ,\mu]$ and $v_\Omega^-[S_{n,\lambda} ,\mu]$ admit a continuous extension to $\overline{\Omega}$ and to $\overline{\Omega^-}$, respectively, we still use the symbols $v_\Omega^+[S_{n,\lambda} ,\mu]$ and $v_\Omega^-[S_{n,\lambda} ,\mu]$ for the continuous extensions and if the values of $v_\Omega^\pm[S_{n,\lambda} ,\mu](x)$ coincide for each $x\in\partial\Omega$, then we set
\[
V_\Omega[S_{n,\lambda} ,\mu](x)\equiv v_\Omega^+[S_{n,\lambda} ,\mu](x)=v_\Omega^-[S_{n,\lambda} ,\mu]^-(x)\qquad\forall x\in\partial\Omega\,.
\]
If $\mu$ is continuous, then it is known that $v_\Omega[S_{n,\lambda} ,\mu ]$ is continuous in ${\mathbb{R}}^n$. Indeed,  $\partial\Omega$ is upper $(n-1)$-Ahlfors regular with respect to ${\mathbb{R}}^n$ and $S_{n,\lambda}$ has a weak singularity (cf. \textit{e.g.}, \cite[Prop.~6.5]{La24}, \cite[Prop.~4.3]{La22b}, \cite[Lem.~4.2 (i)]{DoLa17}). If $\lambda=0$, \textit{i.e.} if $\Delta+\lambda $ is the Laplace operator, and if we take the fundamental solution $S_n$ of $\Delta$ that is delivered by the formula
 \[
S_{n}(\xi)\equiv
\left\{
\begin{array}{lll}
\frac{1}{s_{n}}\ln  |\xi| \qquad &   \forall \xi\in 
{\mathbb{R}}^{n}\setminus\{0\},\quad & {\mathrm{if}}\ n=2\,,
\\
\frac{1}{(2-n)s_{n}}|\xi|^{2-n}\qquad &   \forall \xi\in 
{\mathbb{R}}^{n}\setminus\{0\},\quad & {\mathrm{if}}\ n>2\,,
\end{array}
\right.
\]
where $s_{n}$ denotes the $(n-1)$ dimensional measure of 
$\partial{\mathbb{B}}_{n}(0,1)$, then  we set
\begin{equation}\label{eq:vosn}
v_\Omega[\mu]=v_\Omega[S_n,\mu]\,,\qquad
V_\Omega[\mu]=V_\Omega[S_n,\mu]
\,.
\end{equation}
  Next we introduce the following (classical) technical statement on the  double layer potential.
\begin{theorem}\label{thm:dlay}
 Let  $\alpha\in]0,1[$. Let $\Omega$ be a bounded open subset of  ${\mathbb{R}}^n$ of class $C^{1,\alpha}$. Let $\lambda\in {\mathbb{C}}$.  Let $S_{n,\lambda} $ be a fundamental solution 
of the operator $\Delta+\lambda$. If $\mu\in C^{0,\alpha}(\partial\Omega)$, and
\begin{equation}\label{thm:dlay1}
w_\Omega[S_{n,\lambda} ,\mu](x)\equiv\int_{\partial\Omega}\frac{\partial}{\partial\nu_{\Omega,y}}\left(S_{n,\lambda} (x-y)\right)\mu(y)\,d\sigma_y\qquad\forall x\in {\mathbb{R}}^n\,,
\end{equation}
where
\[
\frac{\partial}{\partial \nu_{\Omega,y} }
\left(S_{n,\lambda}(x-y)\right)\equiv
  -DS_{n,\lambda}(x-y) \nu_{\Omega}(y) \qquad\forall (x,y)\in\mathbb{R}^n\times\partial\Omega\,, x\neq y\,,
 \]
then the restriction 
$w_\Omega[S_{n,\lambda} ,\mu]_{|\Omega}$ can be extended uniquely to a function  
$w^{+}_\Omega [S_{n,\lambda} ,\mu]$ of class $ C^{0,\alpha}(\overline{\Omega})$ and $w_\Omega[S_{n,\lambda} ,\mu]_{|{\mathbb{R}}^n\setminus\overline{\Omega}}$ can be extended uniquely to a function  
$w^{-}_\Omega [S_{n,\lambda} ,\mu]$ of class  $
C^{0,\alpha}_{ {\mathrm{loc}} }(\overline{\Omega^{-}})$.\par
Moreover, the map 
from   $C^{0,\alpha}(\partial\Omega)$ to    $C^{0,\alpha}(\overline{\Omega})$
 which takes $\mu$ to $
w^{+}_\Omega [S_{n,\lambda} ,\mu]$ is  continuous and the map from     $C^{0,\alpha}(\partial\Omega)$ to   $C^{0,\alpha}(\overline{{\mathbb{B}}_{n}(0,r)}\setminus \Omega) $ which takes $\mu$ to 
$w^{-}_\Omega [S_{n,\lambda} ,\mu]_{|\overline{{\mathbb{B}}_{n}(0,r)}\setminus \Omega}$ is continuous for all $r\in]0,+\infty[$ such that $\overline{\Omega}\subseteq {\mathbb{B}}_{n}(0,r)$ and we have the following jump relation
\begin{equation}\label{thm:dlay1a}
w^{\pm}_\Omega [S_{n,\lambda},\mu](x)
=\pm\frac{1}{2}\mu(x)+w_\Omega[S_{n,\lambda},\mu](x)
\qquad\forall x\in\partial\Omega\,.
\end{equation}
\end{theorem}
{\bf Proof.} By definition of double layer potential, we have
\begin{equation}\label{thm:dlay2}
w_{\Omega}[S_{n,\lambda} ,\mu](x)  =-\sum_{j=1}^{n}\frac{\partial}{\partial x_{j}}v_{\Omega}
[S_{n,\lambda} ,\mu  \nu_{j}]  (x)\qquad\forall x\in  {\mathbb{R}}^n\setminus\partial\Omega\,,
\end{equation}
for all $\mu\in C^{0,\alpha}(\partial\Omega)$. 
Since the components of $\nu$ belong to $C^{0,\alpha}(\partial\Omega)$, the continuity of the pointwise product 
from 
 $C^{0,\alpha}(\partial\Omega)\times C^{0,\alpha}(\partial\Omega)$ to $C^{0,\alpha}(\partial\Omega)$
 and the classical Schauder  regularity and jump properties of the single layer potential imply the validity of the statements for $w^{\pm}_\Omega [S_{n,\lambda},\mu]$ (cf.~\textit{e.g.}, \cite[Thm.~7.1]{DoLa17}). \hfill  $\Box$ 

\vspace{\baselineskip}

We also set
\[
W_\Omega[S_{n,\lambda} ,\mu](x)\equiv w_\Omega[S_{n,\lambda} ,\mu](x)\qquad\forall x\in\partial\Omega\,.
\]
If $\lambda=0$, \textit{i.e.} if $\Delta+\lambda $ is the Laplace operator, and if we take the fundamental solution $S_n$ of $\Delta$ , then  we set
\begin{equation}\label{eq:wosn}
w_\Omega[\mu]\equiv w_\Omega[S_n,\mu]\,,\qquad W_\Omega[\mu]\equiv W_\Omega[S_n,\mu]\,.
\end{equation}
\begin{theorem}
\label{thm:rfdlap}
Let  $\alpha\in]0,1[$.  Let $\Omega$ be a bounded open subset of ${\mathbb{R}}^{n}$ of class $C^{1,\alpha}$. Let  $\lambda\in {\mathbb{C}}$. Let $S_{n,\lambda} $ be a fundamental solution 
of the operator $\Delta+\lambda$. If $\mu\in C^{0,\alpha}(\partial\Omega)$, then
\begin{eqnarray}\label{thm:rfdlap1}
\lefteqn{
{\mathcal{G}}_{d,+}[\mu](x)=\lambda\int_{\Omega}{\mathcal{G}}_{d,+}[\mu](y)S_{n,\lambda} (x-y)\,dy
}
\\ \nonumber
&&\qquad
+\int_{\partial\Omega}\mu(y)\frac{\partial}{\partial\nu_{\Omega,y}}\left(S_{n,\lambda} (x-y)\right)\,d\sigma_y
\\ \nonumber
&&\qquad
-\langle r_{|\partial\Omega}^tS_{\Omega,+}^t[\mu](y),S_{n,\lambda} (x-y)\rangle \qquad\forall x\in\Omega\,,
\\ \label{thm:rfdlap1a}
\lefteqn{
0=\lambda\int_{\Omega}{\mathcal{G}}_{d,+}[\mu](y)S_{n,\lambda} (x-y)\,dy
}
\\ \nonumber
&&\qquad
+\int_{\partial\Omega}\mu(y)\frac{\partial}{\partial\nu_{\Omega,y}}\left(S_{n,\lambda} (x-y)\right)\,d\sigma_y
\\ \nonumber
&&\qquad
-\langle r_{|\partial\Omega}^tS_{\Omega,+}^t[\mu](y),S_{n,\lambda} (x-y)\rangle \qquad\forall x\in\Omega^-\,.
\end{eqnarray}
\end{theorem}
{\bf Proof.} By a known approximation result, there exists a sequence
$\{\mu_j\}_{j\in {\mathbb{N}}}$ in $C^{1,\alpha}(\partial\Omega)$ such that
\[
\sup_{j\in {\mathbb{N}}}\|\mu_j\|_{C^{0,\alpha}(\partial\Omega)}<+\infty\,,
\qquad
\lim_{j\to\infty}\mu_j=\mu\quad\text{in}\ C^{0,\beta}(\partial\Omega)\quad\forall \beta\in]0,\alpha[\,,
\]
(cf.~\textit{e.g.}, \cite[Lem.~A.25]{La24b}). Then classical properties of the (harmonic) Green operator imply that ${\mathcal{G}}_{d,+}[\mu_j]\in C^{1,\alpha}(\overline{\Omega})$
(cf.~Theorem \ref{thm:idwp})
 and accordingly  the classical third Green Identity for the Helmholtz operator of Theorem \ref{thm:thirdgreenh} of the Appendix ensures that 
\begin{eqnarray}\label{thm:rfdlap2}
\lefteqn{
{\mathcal{G}}_{d,+}[\mu_j](x)=\lambda\int_{\Omega}{\mathcal{G}}_{d,+}[\mu_j](y)S_{n,\lambda} (x-y)\,dy
}
\\ \nonumber
&&\qquad
+\int_{\partial\Omega}\mu_j(y)\frac{\partial}{\partial\nu_{\Omega,y}}\left(S_{n,\lambda} (x-y)\right)\,dy
\\ \nonumber
&&\qquad
-\int_{\partial\Omega}\frac{\partial}{\partial\nu_\Omega}{\mathcal{G}}_{d,+}[\mu_j](y) S_{n,\lambda} (x-y)\,d\sigma_y\quad\forall x\in\Omega 
\end{eqnarray}
and that
\begin{eqnarray}\label{thm:rfdlap2a}
\lefteqn{
0=\lambda\int_{\Omega}{\mathcal{G}}_{d,+}[\mu_j](y)S_{n,\lambda} (x-y)\,dy
}
\\ \nonumber
&&\qquad
+\int_{\partial\Omega}\mu_j(y)\frac{\partial}{\partial\nu_{\Omega,y}}\left(S_{n,\lambda} (x-y)\right)\,dy
\\ \nonumber
&&\qquad
-\int_{\partial\Omega}\frac{\partial}{\partial\nu_\Omega}{\mathcal{G}}_{d,+}[\mu_j](y) S_{n,\lambda} (x-y)\,d\sigma_y\quad\forall x\in\Omega^-\,,
\end{eqnarray}
for all $j\in {\mathbb{N}}$. By known properties of the (harmonic) Green operator,  we have
\begin{equation}\label{thm:rfdlap2b}
\lim_{j\to\infty}{\mathcal{G}}_{d,+}[\mu_j]={\mathcal{G}}_{d,+}[\mu] 
\qquad\text{in}\  C^{0,\beta}(\overline{\Omega})  \quad\forall \beta\in]0,\alpha[
\end{equation}
(cf.~Theorem \ref{thm:idwp}). Then Theorem \ref{thm:nwtdma}
on the  acoustic volume potential associated to $S_{n,\lambda} $ implies that 
\begin{eqnarray}\nonumber
\lefteqn{
\lim_{j\to\infty}
\int_{\Omega}{\mathcal{G}}_{d,+}[\mu_j](y)S_{n,\lambda} (\cdot-y)\,dy
}
\\ \nonumber
&&\qquad 
=
\int_{\Omega}{\mathcal{G}}_{d,+}[\mu](y)S_{n,\lambda} (\cdot-y)\,dy \ \text{in}\ C^{2,\beta}(\overline{\Omega})\,,
\\ \nonumber
\lefteqn{
\lim_{j\to\infty}
\int_{\Omega}{\mathcal{G}}_{d,+}[\mu_j](y)S_{n,\lambda} (\cdot-y)\,dy 
}
\\  \label{thm:rfdlap3}
&&\qquad 
=
\int_{\Omega}{\mathcal{G}}_{d,+}[\mu](y)S_{n,\lambda} (\cdot-y)\,dy\ \text{in}\ C^{2,\beta}(\overline{{\mathbb{B}}_n(0,r)}\setminus\Omega) \,,
\end{eqnarray}
  for all $r\in]0,+\infty[$ such that $\overline{\Omega}\subseteq {\mathbb{B}}_n(0,r)$ 
and  Theorem \ref{thm:dlay} on the double layer potential implies that 
\begin{equation}\label{thm:rfdlap4}
\lim_{j\to\infty}\int_{\partial\Omega}\mu_j(y)\frac{\partial}{\partial\nu_{\Omega,y}}\left(S_{n,\lambda} (\cdot-y)\right)\,dy
=
\int_{\partial\Omega}\mu(y)\frac{\partial}{\partial\nu_{\Omega,y}}\left(S_{n,\lambda}  (\cdot-y)\right)\,dy
\end{equation}
both in $C^{0,\beta}(\overline{\Omega})$ and in $C^{0,\beta}(\overline{{\mathbb{B}}_n(0,r)}\setminus\Omega)$ for all $r\in]0,+\infty[$ such that $\overline{\Omega}\subseteq {\mathbb{B}}_n(0,r)$ and for all $\beta\in]0,\alpha[$. Next we note that the continuity of $S_{\Omega,+}^t$ from 
$C^{0,\beta}(\partial\Omega)$ to $\left(C^{1,\beta}(\partial\Omega)\right)'$
  implies that 
\[
\lim_{j\to\infty}S_{\Omega,+}^t[\mu_j]=S_{\Omega,+}^t[\mu]\qquad\text{in}\ \left(C^{1,\beta}(\partial\Omega)\right)' \,,
\]
 and accordingly that
\begin{eqnarray}\label{thm:rfdlap5}
\lefteqn{
\lim_{j\to\infty}
\int_{\partial\Omega}\frac{\partial}{\partial\nu_\Omega}{\mathcal{G}}_{d,+}[\mu_j](y) S_{n,\lambda} (x-y)\,d\sigma_y
}
\\ \nonumber
&&
=\lim_{j\to\infty}\langle S_{\Omega,+}^t [\mu_j](y),r_{|\partial\Omega}(S_{n,\lambda} (x-y))\rangle 
\\ \nonumber
&&
=\langle S_{\Omega,+}^t [\mu](y),r_{|\partial\Omega}(S_{n,\lambda} (x-y))\rangle 
=\langle r_{|\partial\Omega}^tS_{\Omega,+}^t [\mu](y),S_{n,\lambda} (x-y)\rangle 
\end{eqnarray}
for all $x\in {\mathbb{R}}^n\setminus\partial\Omega$ and $\beta\in]0,\alpha[$.
 Then by the limiting relations of 
(\ref{thm:rfdlap2b})--(\ref{thm:rfdlap5}), we can take the limit as $j$ tends to infinity in equality (\ref{thm:rfdlap2}) and obtain equality (\ref{thm:rfdlap1})
and in equality (\ref{thm:rfdlap2a}) and obtain equality (\ref{thm:rfdlap1a}).\hfill  $\Box$ 

\vspace{\baselineskip}

We are now ready to prove the following continuity statement for the acoustic single layer potential.
\begin{theorem}\label{thm:slhco-1a} 
Let $\alpha\in]0,1[$. Let $\Omega$ be a bounded open subset of ${\mathbb{R}}^n$ of class $C^{1,\alpha}$. Let $\lambda\in {\mathbb{C}}$.  Let $S_{n,\lambda} $ be a fundamental solution 
of the operator $\Delta+\lambda$. Let $r\in]0,+\infty[$ be such that $\overline{\Omega}\subseteq 
{\mathbb{B}}(0,r)$.

 If $\mu\in V^{-1,\alpha}(\partial\Omega)$, then $v_\Omega^+[S_{n,\lambda} ,\mu]\in C^{0,\alpha}(\overline{\Omega})_\Delta$ and $v_\Omega^-[S_{n,\lambda} ,\mu]\in C^{0,\alpha}_{{\mathrm{loc}}	}(\overline{\Omega^-})_\Delta$. Moreover, the map from $V^{-1,\alpha}(\partial\Omega)$ to $C^{0,\alpha}(\overline{\Omega})_\Delta$ that takes $\mu$ to
 $v_\Omega^+[S_{n,\lambda} ,\mu]$ is continuous 
 and
 the map from $V^{-1,\alpha}(\partial\Omega)$ to $C^{0,\alpha}(\overline{{\mathbb{B}}(0,r)}\setminus\Omega)_\Delta$ that takes $\mu$ to
 $v_\Omega^-[S_{n,\lambda} ,\mu]_{|\overline{{\mathbb{B}}(0,r)}\setminus\Omega}$ is continuous.
 \end{theorem}
 {\bf Proof.} By the definition of $V^{-1,\alpha}(\partial\Omega)$, by the Lemma 13.5  of
 \cite{La24b}
 on the continuity of linear maps defined on
 $V^{-1,\alpha}(\partial\Omega)$  and by the equalities
\begin{eqnarray*}
&&\Delta v_\Omega^+[S_{n,\lambda} ,\mu_0+S_{\Omega,+}^t[\mu_1]]=-\lambda v_\Omega^+[S_{n,\lambda} ,\mu_0+S_{\Omega,+}^t[\mu_1]]\ \  \text{in}\ \Omega\,,
\\
 &&\Delta v_\Omega^-[S_{n,\lambda} ,\mu_0+S_{\Omega,+}^t[\mu_1]] =-\lambda v_\Omega^-[S_{n,\lambda} ,\mu_0+S_{\Omega,+}^t[\mu_1]]  \ \  \text{in}\  {\mathbb{B}}(0,r)\setminus\overline{\Omega}\,,
\end{eqnarray*}
for all $\mu_0$, $\mu_1\in C^{0,\alpha}(\partial\Omega)$, 
 it suffices to show that if $\mu_0$, $\mu_1\in C^{0,\alpha}(\partial\Omega)$, then
\begin{eqnarray}\label{thm:slhco-1a1}
&&v_\Omega^+[S_{n,\lambda} ,\mu_0+S_{\Omega,+}^t[\mu_1]]\in C^{0,\alpha}(\overline{\Omega})\,,\quad
\\\nonumber
&&v_\Omega^-[S_{n,\lambda} ,\mu_0+S_{\Omega,+}^t[\mu_1]]_{|\overline{{\mathbb{B}}(0,r)}\setminus\Omega}
\in C^{0,\alpha}(\overline{{\mathbb{B}}(0,r)}\setminus\Omega)
\end{eqnarray}
and the maps from $\left(C^{0,\alpha}(\partial\Omega)\right)^2$ to $C^{0,\alpha}(\overline{\Omega})$ and to $C^{0,\alpha}(\overline{{\mathbb{B}}(0,r)}\setminus\Omega)$ that take  $(\mu_0,\mu_1)$ to $v_\Omega^+[S_{n,\lambda} ,\mu_0+S_{\Omega,+}^t[\mu_1]]$ and to $v_\Omega^-[S_{n,\lambda} ,\mu_0+S_{\Omega,+}^t[\mu_1]]_{|\overline{{\mathbb{B}}(0,r)}\setminus\Omega}$ are continuous, respectively.

By a classical result on the single layer potential, we know that  
\begin{eqnarray}\label{thm:slhco-1a2}
&&v_\Omega^+[S_{n,\lambda} ,\cdot]\in {\mathcal{L}}\left(C^{0,\alpha}(\partial\Omega),C^{1,\alpha}(\overline{\Omega})\right)\,,
\\	\nonumber
&&v_\Omega^-[S_{n,\lambda} ,\cdot]_{|\overline{{\mathbb{B}}(0,r)}\setminus\Omega}
\in {\mathcal{L}}\left(C^{0,\alpha}(\partial\Omega),C^{1,\alpha}(\overline{{\mathbb{B}}(0,r)}\setminus\Omega)\right)\,,
\end{eqnarray}
(cf.~\textit{e.g.}, \cite[Thm.~7.1]{DoLa17}). Next we note that Theorem \ref{thm:rfdlap} implies that
\begin{eqnarray} \nonumber
\lefteqn{
v_\Omega^+[S_{n,\lambda} ,S_{\Omega,+}^t[\mu_1]](x)=-{\mathcal{G}}_{d,+}[\mu_1](x)+\lambda\int_{\Omega}{\mathcal{G}}_{d,+}[\mu_1](y)S_{n,\lambda} (x-y)\,dy
}
\\ \nonumber
&&\qquad\qquad\qquad
+\int_{\partial\Omega}\mu_1(y)\frac{\partial}{\partial\nu_{\Omega,y}}S_{n,\lambda} (x-y)\,d\sigma_y
 \qquad\forall x\in\Omega\,,
\\ \label{thm:slhco-1a3} 
\lefteqn{
v_\Omega^-[S_{n,\lambda} ,S_{\Omega,+}^t[\mu_1](x)=\lambda\int_{\Omega}{\mathcal{G}}_{d,+}[\mu_1](y)S_{n,\lambda} (x-y)\,dy
}
\\ \nonumber
&&\qquad\qquad\qquad
+\int_{\partial\Omega}\mu_1(y)\frac{\partial}{\partial\nu_{\Omega,y}}S_{n,\lambda} (x-y)\,d\sigma_y
\qquad\forall x\in\Omega^-\,.
\end{eqnarray}
By a classical property of the (harmonic) Green operator, we have
\begin{equation}\label{thm:slhco-1a4}
{\mathcal{G}}_{d,+}[\cdot]\in 
{\mathcal{L}}\left(C^{0,\alpha}(\partial\Omega),C^{0,\alpha}(\overline{\Omega})\right)
\end{equation}
(cf.~Theorem \ref{thm:idwp}). By a classical property of the acoustic volume potential, we have
\begin{eqnarray}\label{thm:slhco-1a5}
&&{\mathcal{P}}_\Omega^+[S_{n,\lambda} ,\cdot]\in 
{\mathcal{L}}\left(C^{0,\alpha}(\overline{\Omega}),C^{2,\alpha}(\overline{\Omega})\right)\,,
\\	\nonumber
&&{\mathcal{P}}_\Omega^-[S_{n,\lambda} ,\cdot]_{|\overline{{\mathbb{B}}(0,r)}\setminus\Omega}\in 
{\mathcal{L}}\left(C^{0,\alpha}(\overline{\Omega}),C^{2,\alpha}(\overline{{\mathbb{B}}(0,r)}\setminus\Omega)\right)\,,
\end{eqnarray}
(cf.~Theorem \ref{thm:nwtdma}).
  Then classical properties of the acoustic double layer potential imply that
\begin{eqnarray}\label{thm:slhco-1a6}
&&w_\Omega^+[S_{n,\lambda} ,\cdot]\in 
{\mathcal{L}}\left(C^{0,\alpha}(\partial\Omega),C^{0,\alpha}(\overline{\Omega})\right)\,,
\\	\nonumber
&&w_\Omega^-[S_{n,\lambda} ,\cdot]_{|\overline{{\mathbb{B}}(0,r)}\setminus\Omega}\in 
{\mathcal{L}}\left(C^{0,\alpha}(\partial\Omega),C^{0,\alpha}(\overline{{\mathbb{B}}(0,r)}\setminus\Omega)\right)\,,
\end{eqnarray}
(cf.~Theorem \ref{thm:dlay}). 
Then the formulas in  (\ref{thm:slhco-1a3}) and the continuity properties of (\ref{thm:slhco-1a4})--(\ref{thm:slhco-1a6}) imply that
\begin{eqnarray}\label{thm:slhco-1a7}
&&v_\Omega^+[S_{n,\lambda} ,S_{\Omega,+}^t[\cdot]]\in 
{\mathcal{L}}\left(C^{0,\alpha}(\partial\Omega),C^{0,\alpha}(\overline{\Omega})\right)\,,
\\	\nonumber
&&v_\Omega^-[S_{n,\lambda} ,S_{\Omega,+}^t[\cdot]]_{|\overline{{\mathbb{B}}(0,r)}\setminus\Omega}\in 
{\mathcal{L}}\left(C^{0,\alpha}(\partial\Omega),C^{0,\alpha}(\overline{{\mathbb{B}}(0,r)}\setminus\Omega)\right)\,.
\end{eqnarray}
Then the continuity properties of (\ref{thm:slhco-1a2}) and (\ref{thm:slhco-1a7}) imply the validity of (\ref{thm:slhco-1a1}) and of the subsequent continuity properties. Hence,  the proof is complete.\hfill  $\Box$ 

\vspace{\baselineskip}

\section{Jump formulas for the  acoustic single layer potential}
\label{sec:silahjump}

We first prove the following elementary embedding lemma.
\begin{lemma}\label{lem:am0b-1a}
 Let $\alpha\in]0,1[$.  Let $\Omega$ be a bounded open subset of ${\mathbb{R}}^n$ of class $C^{1,\alpha}$. If $\beta\in]0,1[$, then $C^{0,\beta}(\partial\Omega)$ is continuously embedded into $V^{-1,\alpha}(\partial\Omega)$.
 \end{lemma}
 {\bf Proof.} If $\beta\geq \alpha$, then  $C^{0,\beta}(\partial\Omega)$ is continuously embedded into $C^{0,\alpha}(\partial\Omega)$ and  Definition \ref{defn:v-1a} ensures that 
$C^{0,\alpha}(\partial\Omega)$  is continuously embedded into $V^{-1,\alpha}(\partial\Omega)$. 
We now consider case $\beta<\alpha$ and we first prove that $C^{0,\beta}(\partial\Omega)$ is contained in $V^{-1,\alpha}(\partial\Omega)$. If $\tau\in C^{0,\beta}(\partial\Omega)$, then a classical  regularity result on the harmonic single layer potential implies that
\begin{equation}\label{lem:am0b-1a1}
v_\Omega^+[\tau]\in C^{1,\beta}(\overline{\Omega})\subseteq C^{0,\alpha}(\overline{\Omega})\,,\ \ 
v_\Omega^-[\tau]\in C^{1,\beta}_{	{\mathrm{loc}}	}(\overline{\Omega^-})\subseteq C^{0,\alpha}_{	{\mathrm{loc}}	}(\overline{\Omega^-})\, \ \ v_\Omega^+[\tau]_{|\partial\Omega}=v_\Omega^-[\tau]_{|\partial\Omega}
\end{equation}
(cf.~\textit{e.g.}, \cite[Thm.~7.1]{DoLa17})
and the Fubini Theorem implies that
\begin{equation}\label{lem:am0b-1a2}
\int_{\partial\Omega}\tau v_\Omega^+[\psi]\,d\sigma
=\int_{\partial\Omega}\psi v_\Omega^+[\tau]\,d\sigma\qquad\forall\psi\in C^{0,\alpha}(\partial\Omega)\,.
\end{equation}
Then (\ref{lem:am0b-1a1}), (\ref{lem:am0b-1a2}) and a   result of \cite[Thm.~18.1 (see also Defns.~13.2, 15.10)]{La24b} on $V^{-1,\alpha}(\partial\Omega)$ imply that $\tau\in 
V^{-1,\alpha}(\partial\Omega)$. Hence, we have proved that $C^{0,\beta}(\partial\Omega)$ is contained in $ V^{-1,\alpha}(\partial\Omega)$.

 By a criterium of \cite[Lem.~18.7]{La24b},  the inclusion of $C^{0,\beta}(\partial\Omega)$ into $ V^{-1,\alpha}(\partial\Omega)$ is continuous if and only if the map from 
$C^{0,\beta}(\partial\Omega)$ to $C^{0,\alpha}(\partial\Omega)$ that takes $\tau$ to 
\begin{equation}\label{lem:am0b-1a3}
v_\Omega^+\left[\tau-\frac{\int_{\partial\Omega}\tau\,d\sigma}{\int_{\partial\Omega} \,d\sigma}\right]_{|\partial\Omega}+\frac{\int_{\partial\Omega}\tau\,d\sigma}{\int_{\partial\Omega} \,d\sigma}
\end{equation}
is continuous. By a classical result, the harmonic single layer potential $v_\Omega^+[\cdot]_{|\partial\Omega}$ is continuous from $C^{0,\beta}(\partial\Omega)$ to $C^{1,\beta}(\partial\Omega)$ that is continuously embedded into $C^{0,\alpha}(\partial\Omega)$. 
Since the map from $C^{0,\beta}(\partial\Omega)$  to  both
$C^{0,\beta}(\partial\Omega)$ and $C^{0,\alpha}(\partial\Omega)$
that takes $\tau$ to the constant  $\frac{\int_{\partial\Omega}\tau\,d\sigma}{\int_{\partial\Omega} \,d\sigma}$
is continuous, then map in (\ref{lem:am0b-1a3}) is continuous from  $C^{0,\beta}(\partial\Omega)$  to  $C^{0,\alpha}(\partial\Omega)$ and thus the proof is complete.\hfill  $\Box$ 

\vspace{\baselineskip}

Next we introduce an approximation lemma for the distributions in $V^{-1,\alpha}(\partial\Omega)$. 
\begin{lemma}\label{lem:apr-1a}
 Let $\alpha\in]0,1[$. 
 Let $\Omega$ be a bounded open subset of ${\mathbb{R}}^n$ of class $C^{1,\alpha}$. 
If $\tau\in V^{-1,\alpha}(\partial\Omega)$, then there exists a sequence $\{\tau_j\}_{j\in {\mathbb{N}}}$ in 
 $C^{1,\alpha}(\partial\Omega)$ such that
 \begin{equation}\label{lem:apr-1a1}
 \sup_{j\in {\mathbb{N}}}\|\tau_j\|_{
 V^{-1,\alpha}(\partial\Omega) 
 }<+\infty\,,\qquad
 \lim_{j\to\infty}\tau_j=\tau\quad\text{in}\ V^{-1,\beta}(\partial\Omega)\quad \forall\beta\in]0,\alpha[\,.
 \end{equation}
\end{lemma}
{\bf Proof.} Let $\mu_0$, $\mu_1\in C^{0,\alpha}(\partial\Omega)$ be such that
\begin{equation}
\tau=\mu_0+S_{\Omega,+}^t[\mu_1]\,.
\end{equation}
A known approximation property implies that there exists a sequence $\{\mu_{s,l}\}_{l\in {\mathbb{N}}}$ in $C^{1,\alpha}(\partial\Omega)$   that converges to $\mu_s$ in the $ C^{0,\beta}(\partial\Omega)$-norm for all $\beta\in]0,\alpha[$ and that is bounded in the $ C^{0,\alpha}(\partial\Omega)$-norm, for each $s\in\{0,1\}$   (cf.~\textit{e.g.}, \cite[Lem.~A.25]{La24b}). Since the map $A$ from $\left(C^{0,\beta}(\partial\Omega)\right)^2$ to $V^{-1,\beta}(\partial\Omega)$ that takes a pair $(\theta_0,\theta_1)$ to $\theta_0+S_{\Omega,+}^t[\theta_1]$ is continuous (cf.~Definition \ref{defn:v-1a}), we have
 \begin{equation}\label{lem:apr1a2}
\lim_{l\to\infty}\mu_{0,l}+S_{\Omega,+}^t[\mu_{1,l}]=\mu_{0}+S_{\Omega,+}^t[\mu_{1}]=\tau\qquad
\text{in}\quad V^{-1,\beta}(\partial\Omega)\,.
\end{equation}
Since $A$ is linear and continuous from $\left(C^{0,\alpha}(\partial\Omega)\right)^2$ to $V^{-1,\alpha}(\partial\Omega)$, we also have
\begin{equation}\label{lem:apr1a3}
\sup_{l\in{\mathbb{N}} }\|\mu_{0,l}+S_{\Omega,+}^t[\mu_{1,l}]\|_{V^{-1,\alpha}(\partial\Omega)} <+\infty\,.
\end{equation}
By Theorem \ref{thm:idwp} and \cite[Thm.~5.9 (ii)]{La24b}, we have 
\[
S_{\Omega,+}^t[\theta]=\frac{\partial}{\partial\nu}{\mathcal{G}}_{d,+}[\theta]\in C^{0,\alpha}(\partial\Omega)
\qquad\forall\theta\in C^{1,\alpha}(\partial\Omega)
\]
and  $S_{\Omega,+}^t$ is linear and continuous from $C^{1,\alpha}(\partial\Omega)$ to 
$C^{0,\alpha}(\partial\Omega)$. Then we set $\sigma_l\equiv \mu_{0,l}+S_{\Omega,+}^t[\mu_{1,l}]$ for all $l\in {\mathbb{N}}$ and we note that $\sigma_l\in C^{0,\alpha}(\partial\Omega)$ for all $l\in {\mathbb{N}}$ and that
\begin{equation}\label{lem:apr1a3a}
 \sup_{l\in {\mathbb{N}}}\|\sigma_l\|_{
 V^{-1,\alpha}(\partial\Omega) 
 }<+\infty\,,\qquad
 \lim_{l\to\infty}\sigma_l=\tau\quad\text{in}\ V^{-1,\beta}(\partial\Omega)\quad \forall\beta\in]0,\alpha[\,.
 \end{equation}
 Here we note that  $\sigma_l\in C^{0,\alpha}(\partial\Omega)$ for all $l\in {\mathbb{N}}$, but we need a sequence $\{\tau_j\}_{j\in {\mathbb{N}}}$ in $C^{1,\alpha}(\partial\Omega)$ that satisfies the conditions of (\ref{lem:apr-1a1}). In order to prove the existence of $\{\tau_j\}_{j\in {\mathbb{N}}}$,
 we now fix $\tilde{\beta}\in ]0,\alpha[$ and  define a sequence $\{\tau_j\}_{j\in {\mathbb{N}}}$  in $C^{1,\alpha}(\partial\Omega)$ that satisfies the conditions in 
 (\ref{lem:apr-1a1}) for $\beta=\tilde{\beta}$ and then show that the same sequence satisfies the conditions in 
 (\ref{lem:apr-1a1}) for all $\beta\in]0,\alpha[$.  Let $j\in {\mathbb{N}}$. By 
(\ref{lem:apr1a3a}) with  $\beta=\tilde{\beta}$, there exists $l_j\in{\mathbb{N}}$ such that
\begin{equation}\label{lem:apr1a4}
\|\tau-\sigma_{l_j}\|_{ V^{-1,\tilde{\beta}}(\partial\Omega) }\leq 2^{-j-1}\,.
 \end{equation}
Since $\sigma_{l_j}\in C^{0,\alpha}(\partial\Omega)$, a known approximation lemma implies that there exists $\tau_j\in C^{1,\alpha}(\partial\Omega)$ such that
\begin{equation}\label{lem:apr1a5}
\|\sigma_{l_j}-\tau_j\|_{ V^{-1,\tilde{\beta}}(\partial\Omega) }
\leq
\|J\|_{
{\mathcal{L}}\left(C^{0,\tilde{\beta}}(\partial\Omega), V^{-1,\tilde{\beta}}(\partial\Omega) \right)}
\|\sigma_{l_j}-\tau_j\|_{ C^{0,\tilde{\beta}}(\partial\Omega) }\leq 2^{-j-1}\,,
 \end{equation}
 where $J$ denotes the inclusion map from $C^{0,\tilde{\beta}}(\partial\Omega)$ to $ V^{-1,\tilde{\beta}}(\partial\Omega) $ (cf.~\cite[Lem.~A.25]{La24b}, Lemma \ref{lem:am0b-1a}). By (\ref{lem:apr1a4}) and  (\ref{lem:apr1a5}), we have
\begin{eqnarray}
\lefteqn{\|\tau-\tau_j\|_{ V^{-1,\tilde{\beta}}(\partial\Omega) }\leq 2^{-j}\,,
}
\\ \nonumber
\lefteqn{
\|\tau_j\|_{ V^{-1,\alpha}(\partial\Omega) }
\leq \|\tau_j-\sigma_{l_j}\|_{ V^{-1,\alpha}(\partial\Omega) }
+\|\sigma_{l_j}\|_{ V^{-1,\alpha}(\partial\Omega) }
}
\\ \nonumber
&&\qquad
\leq
\|J\|_{
{\mathcal{L}}\left(C^{0,\tilde{\beta}}(\partial\Omega), V^{-1,\alpha}(\partial\Omega) \right)}\|\tau_j-\sigma_{l_j}\|_{C^{0,\tilde{\beta}}(\partial\Omega)}
+ \sup_{l\in {\mathbb{N}}}\|\sigma_l\|_{
 V^{-1,\alpha}(\partial\Omega) 
 }
\\ \nonumber
&&\qquad
\leq
\|J\|_{
{\mathcal{L}}\left(C^{0,\tilde{\beta}}(\partial\Omega), V^{-1,\alpha}(\partial\Omega) \right)}+ \sup_{l\in {\mathbb{N}}}\|\sigma_l\|_{
 V^{-1,\alpha}(\partial\Omega)} <+\infty\,,
\end{eqnarray}
for all $j\in {\mathbb{N}}$, where $J$ denotes the inclusion map from $C^{0,\tilde{\beta}}(\partial\Omega)$ to $ V^{-1,\alpha}(\partial\Omega) $ (cf.~Lemma \ref{lem:am0b-1a}). Hence, the sequence $\{\tau_j\}_{j\in {\mathbb{N}}}$  satisfies the conditions in 
 (\ref{lem:apr-1a1}) for $\beta=\tilde{\beta}$. Next we fix an arbitrary $\beta\in]0,\alpha[$ prove that $\{\tau_j\}_{j\in {\mathbb{N}}}$ converges to $\tau$ also in the $ V^{-1,\beta}(\partial\Omega) $-norm. By a known contradiction argument, it suffices to show that if $\{\tau_{j_k}\}_{k\in {\mathbb{N}}}$ is an arbitrary subsequence of $\{\tau_j\}_{j\in {\mathbb{N}}}$, then 
$\{\tau_{j_k}\}_{k\in {\mathbb{N}}}$ has a subsequence that converges to $\tau$ in the $ V^{-1,\beta}(\partial\Omega) $-norm. Since $ V^{-1,\alpha}(\partial\Omega) $ is compactly embedded into
$ V^{-1,	\max\{\beta,\tilde{\beta}\}	}(\partial\Omega) $ and $\{\tau_{j_k}\}_{k\in {\mathbb{N}}}$ is bounded in the $ V^{-1,\alpha}(\partial\Omega) $-norm, then $\{\tau_{j_k}\}_{k\in {\mathbb{N}}}$ has a subsequence that converges to some $\xi\in V^{-1,\max\{\beta,\tilde{\beta}\}}(\partial\Omega) $ in the $ V^{-1,\max\{\beta,\tilde{\beta}\}}(\partial\Omega) $-norm (cf.~\cite[Thm.~18.8]{La24b}). Since 
$V^{-1,\max\{\beta,\tilde{\beta}\}}(\partial\Omega) $ is continuously embedded into $V^{-1,\tilde{\beta}}(\partial\Omega) $ (cf.~\cite[Thm.~18.8]{La24b}) and  $\{\tau_{j_k}\}_{k\in {\mathbb{N}}}$ converges to $\tau$ in $V^{-1,\tilde{\beta}}(\partial\Omega) $, we conclude that $\tau=\xi\in V^{-1,\max\{\beta,\tilde{\beta}\}}(\partial\Omega) $. Since $V^{-1,\max\{\beta,\tilde{\beta}\}}(\partial\Omega) $ is continuously embedded into $V^{-1,\beta}(\partial\Omega) $ (cf.~\cite[Thm.~18.8]{La24b}), we conclude that $\{\tau_{j_k}\}_{k\in {\mathbb{N}}}$ has a subsequence that converges to $\tau$ in the $ V^{-1,\beta}(\partial\Omega) $-norm. Hence, the sequence $\{\tau_j\}_{j\in {\mathbb{N}}}$  satisfies the conditions in 
 (\ref{lem:apr-1a1}) for $\beta$ too  and thus the proof is complete.\hfill  $\Box$ 

\vspace{\baselineskip}

We are now ready to prove the jump formulas for the normal derivative of the single layer acoustic potential.
\begin{theorem}\label{thm:jufonslap}
 Let $\alpha\in]0,1[$. Let $\Omega$ be a bounded open subset of ${\mathbb{R}}^n$ of class $C^{1,\alpha}$. Let $\lambda\in {\mathbb{C}}$.  Let $S_{n,\lambda} $ be a fundamental solution 
of the operator $\Delta+\lambda$. 

 If $\tau\in V^{-1,\alpha}(\partial\Omega) $, then the following jump formulas hold.
 \begin{eqnarray}
  \label{thm:jufonslap0}
 &&v_\Omega^+[S_{n,\lambda} ,\tau]_{|\partial\Omega}=v_\Omega^-[S_{n,\lambda} ,\tau]_{|\partial\Omega}  
 \\
 \label{thm:jufonslap1}
&&\partial_\nu v_\Omega^+[S_{n,\lambda} ,\tau]=-\frac{1}{2}\tau+W_\Omega^t[S_{n,\lambda} ,\tau]\,,
\\	\label{thm:jufonslap2}
&&-\partial_{\nu_{\Omega^-}} v_\Omega^-[S_{n,\lambda} ,\tau]=+\frac{1}{2}\tau+W_\Omega^t[S_{n,\lambda} ,\tau]\,.
\end{eqnarray}
\end{theorem}
{\bf Proof.}  By Lemma \ref{lem:apr-1a}, there exists a sequence  $\{\tau_j\}_{j\in {\mathbb{N}} }$ as in (\ref{lem:apr-1a1}). Let $\beta\in]0,\alpha[$,  $v\in C^{1,\alpha}(\partial\Omega)$. 
 By the continuity Theorem \ref{thm:slhco-1a}, we have
\begin{eqnarray*}
&&v_\Omega^+[S_{n,\lambda} ,\tau]=\lim_{j\to\infty}v_\Omega^+[S_{n,\lambda} ,\tau_j]   
\qquad\text{in}\ C^{0,\beta}(\overline{\Omega})_\Delta\,,
\\
&&v_\Omega^-[S_{n,\lambda} ,\tau]=\lim_{j\to\infty}v_\Omega^-[S_{n,\lambda} ,\tau_j]   
\qquad\text{in}\ C^{0,\alpha}_{	{\mathrm{loc}}	}(\overline{\Omega^-})_\Delta\,,
\end{eqnarray*}
(see Lemma \ref{lem:c1alcof} for the topology of $C^{0,\alpha}_{
{\mathrm{loc}}	}(\overline{\Omega^-})_\Delta $).   Since $\tau_j\in C^{1,\alpha}(\partial\Omega)\subseteq C^{0,\alpha}(\partial\Omega)$, a classical regularity result on the single layer potential implies that $v_\Omega [S_{n,\lambda} ,\tau_j]$ is continuous in the whole of ${\mathbb{R}}^n$ and that accordingly
$v_\Omega^+[S_{n,\lambda} ,\tau_j]_{|\partial\Omega}=v_\Omega^-[S_{n,\lambda} ,\tau_j]_{|\partial\Omega} $ for all $j\in {\mathbb{N}}$ (cf.~\textit{e.g.}, \cite[Thm~7.1]{DoLa17}).
Then the above limiting relations imply that the validity of formula 
 (\ref{thm:jufonslap0}).\par

 Then the continuity of $\partial_\nu$ from $C^{0,\beta}(\overline{\Omega})_\Delta$ to
$V^{-1,\beta}(\partial\Omega)$ of Proposition \ref{prop:ricodnu}, the continuity of $\partial_{\nu_{\Omega^-}}$ from $C^{0,\alpha}_{
{\mathrm{loc}}	}(\overline{\Omega^-})_\Delta$ to
$V^{-1,\beta}(\partial\Omega)$ of Proposition \ref{prop:recodnu}, 
the inclusion of 
$V^{-1,\beta}(\partial\Omega)$ into the dual of $C^{1,\beta}(\partial\Omega)$
and of $C^{1,\alpha}(\partial\Omega)$ into $C^{1,\beta}(\partial\Omega)$
imply that
\begin{eqnarray}\label{thm:jufonslap3}
&&\langle \partial_\nu v_\Omega^+[S_{n,\lambda} ,\tau],v\rangle=\lim_{j\to\infty}\langle \partial_\nu v_\Omega^+[S_{n,\lambda} ,\tau_j],v\rangle 
\,,
\\ \nonumber
&&
\langle \partial_{\nu_{\Omega^-}} v_\Omega^-[S_{n,\lambda} ,\tau],v\rangle=\lim_{j\to\infty}\langle \partial_{\nu_{\Omega^-}} v_\Omega^-[S_{n,\lambda} ,\tau_j],v\rangle \,.
\end{eqnarray}
Since $\tau_j\in C^{1,\alpha}(\partial\Omega)\subseteq C^{0,\alpha}(\partial\Omega)$, a classical regularity result on the single layer potential implies that $v_\Omega^+[S_{n,\lambda} ,\tau_j]$ belongs to $C^{1,\alpha}(\overline{\Omega})$ for all $j\in {\mathbb{N}}$ (cf.~\textit{e.g.}, \cite[Thm~7.1]{DoLa17}). Then the classical   jump formula  
and the membership    of $W_\Omega[S_{n,\lambda} ,v]$ in $C^{1,\beta}(\partial\Omega)$ 
imply that
\begin{eqnarray}\label{thm:jufonslap6}
\lefteqn{
\lim_{j\to\infty}\langle \partial_\nu v_\Omega^+[S_{n,\lambda} ,\tau_j],v\rangle 
=\lim_{j\to\infty}
\int_{\partial\Omega}
\frac{\partial}{\partial\nu} v_\Omega^+[S_{n,\lambda} ,\tau_j] v\,d\sigma
}
\\ \nonumber
&&\qquad
=\lim_{j\to\infty}
\int_{\partial\Omega}\left(-\frac{1}{2}\tau_j+W_\Omega^t[S_{n,\lambda} ,\tau_j]\right)v\,d\sigma
\\ \nonumber
&&\qquad
=\lim_{j\to\infty}
\langle-\frac{1}{2}\tau_j,v\rangle 
+\lim_{j\to\infty}\langle \tau_j,W_\Omega[S_{n,\lambda} ,v]\rangle
\\ \nonumber
&&\qquad
= 
\langle-\frac{1}{2}\tau,v\rangle 
+ \langle \tau ,W_\Omega[S_{n,\lambda} ,v]\rangle
= 
\langle-\frac{1}{2}\tau +W_\Omega^t[S_{n,\lambda},\tau ] ,  v\rangle
\,,
\end{eqnarray}
(cf.~\textit{e.g.},    \cite[Thm~7.1, Thm.~9.2 (ii)]{DoLa17}). 
Then equalities (\ref{thm:jufonslap3}) and (\ref{thm:jufonslap6}) imply the validity of formula 
 (\ref{thm:jufonslap1}).\par
 
Similarly, we now turn to prove formula (\ref{thm:jufonslap2}). Since $\tau_j\in C^{1,\alpha}(\partial\Omega)\subseteq C^{0,\alpha}(\partial\Omega)$, a classical regularity result on the single layer potential implies that $v_\Omega^-[S_{n,\lambda} ,\tau_j]$ belongs to $C^{1,\alpha}_{{\mathrm{loc}}}(\overline{\Omega^-})$ for all $j\in {\mathbb{N}}$ (cf.~\textit{e.g.}, \cite[Thm~7.1]{DoLa17}). Then the classical   jump formula  
and the membership    of $W_\Omega[S_{n,\lambda} ,v]$ in $C^{1,\beta}(\partial\Omega)$  
imply that
\begin{eqnarray}\label{thm:jufonslap7}
\lefteqn{
\lim_{j\to\infty}\langle \partial_{\nu_{\Omega^-}} v_\Omega^-[S_{n,\lambda} ,\tau_j],v\rangle 
=
\lim_{j\to\infty}\int_{\partial\Omega}\frac{\partial}{\partial\nu_{\Omega^-}}v_\Omega^-[S_{n,\lambda} ,\tau_j]v\,d\sigma
}
\\ \nonumber
&&\qquad
=\lim_{j\to\infty}
\int_{\partial\Omega}\left(-\frac{1}{2}\tau_j-W_\Omega^t[S_{n,\lambda} ,\tau_j]\right)v\,d\sigma
\\ \nonumber
&&\qquad
=\lim_{j\to\infty}
\langle -\frac{1}{2}\tau_j,v\rangle
+\lim_{j\to\infty}\langle \tau_j,-W_\Omega[S_{n,\lambda} ,v]\rangle 
\\ \nonumber
&&\qquad
= 
\langle -\frac{1}{2}\tau,v\rangle
+ \langle \tau ,-W_\Omega[S_{n,\lambda} ,v]\rangle 
= 
\langle -\frac{1}{2}\tau  -W_\Omega^t[S_{n,\lambda} ,\tau] ,  v\rangle
\,,
\end{eqnarray}
(cf.~\textit{e.g.},   \cite[Thm~7.1, Thm.~9.2 (ii)]{DoLa17}). 
 Then equalities  (\ref{thm:jufonslap3}) and	(\ref{thm:jufonslap7}) imply the validity of formula 
 (\ref{thm:jufonslap2}).\hfill  $\Box$ 

\vspace{\baselineskip}

By the jump formula (\ref{thm:jufonslap1})  of Theorem \ref{thm:jufonslap}, by the continuity Theorem \ref{thm:slhco-1a} for $v_\Omega^+[S_{n,\lambda},\cdot]$ and by the continuity of the normal derivative $\partial_{\nu}$ of Proposition \ref{prop:ricodnu}, we deduce the validity of the following statement.
\begin{proposition}\label{prop:cowt-1a}
 Let $\alpha\in]0,1[$. Let $\Omega$ be a bounded open subset of ${\mathbb{R}}^n$ of class $C^{1,\alpha}$. Let $\lambda\in {\mathbb{C}}$.  Let $S_{n,\lambda}$ be a fundamental solution 
of the operator $\Delta+\lambda$. Then the transpose $W_\Omega^t[S_{n,\lambda} ,\cdot]$ of the operator $W_\Omega[S_{n,\lambda} ,\cdot]$ from $C^{1,\alpha}(\partial\Omega)$ to itself  is continuous from  $V^{-1,\alpha}(\partial\Omega) $ to itself. 
\end{proposition}

For the normal derivative of a double layer acoustic potential with H\"{o}lder continuous densities, we have the following statement that generalizes the corresponding known classical statement.
\begin{theorem}\label{thm:junodedl}
Let $\alpha\in]0,1[$. Let $\Omega$ be a bounded open subset of ${\mathbb{R}}^n$ of class $C^{1,\alpha}$. Let $\lambda\in {\mathbb{C}}$.  Let $S_{n,\lambda}$ be a fundamental solution 
of the operator $\Delta+\lambda$. If $\tau\in C^{0,\alpha}(\partial\Omega)$, then the following jump formula holds true
\begin{equation}\label{thm:junodedl1}
\partial_\nu w_\Omega^+[S_{n,\lambda} ,\tau]=-\partial_{\nu_{\Omega^-} }w_\Omega^-[S_{n,\lambda} ,\tau]\,.
\end{equation}
\end{theorem}
{\bf Proof.} By a known approximation result, there exists a sequence
$\{\tau_j\}_{j\in {\mathbb{N}}}$ in $C^{1,\alpha}(\partial\Omega)$ such that
\[
\sup_{j\in {\mathbb{N}}}\|\tau_j\|_{C^{0,\alpha}(\partial\Omega)}<+\infty\,,
\qquad
\lim_{j\to\infty}\tau_j=\mu\quad\text{in}\ C^{0,\beta}(\partial\Omega)\quad\forall \beta\in]0,\alpha[\,,
\]
(cf.~\textit{e.g.}, \cite[Lem.~A.25]{La24b}). Then Theorem \ref{thm:dlay} on the double layer potential and  the equalities
\[
\Delta w_\Omega^\pm[S_{n,\lambda} ,\tau_j]=-\lambda w_\Omega^\pm[S_{n,\lambda} ,\tau_j]\qquad\forall j\in {\mathbb{N}}\,,
\]
imply that
\begin{eqnarray*}
&&\lim_{j\to\infty}w_\Omega^+[S_{n,\lambda} ,\tau_j]=w_\Omega^+[S_{n,\lambda} ,\tau]
\qquad\text{in}\ C^{0,\beta}(\overline{\Omega})_\Delta\,,
\\ \nonumber
&& 
\lim_{j\to\infty}w_\Omega^-[S_{n,\lambda} ,\tau_j]=w_\Omega^-[S_{n,\lambda} ,\tau]
\qquad\text{in}\ C^{0,\beta}_{\mathrm{loc}}(\overline{\Omega^-})_\Delta\,.
\end{eqnarray*}
Then the continuity of $\partial_\nu$ and $\partial_{\nu_{\Omega^-} }$ of Proposition  \ref{prop:ricodnu} and Proposition \ref{prop:recodnu} imply that
\begin{eqnarray}\label{thm:junodedl2}
&&\lim_{j\to\infty}\partial_\nu w_\Omega^+[S_{n,\lambda} ,\tau_j]=\partial_\nu w_\Omega^+[S_{n,\lambda} ,\tau]
\qquad\text{in}\ V^{-1,\beta}(\partial\Omega)\,,
\\ \nonumber
&& 
\lim_{j\to\infty}\partial_{\nu_{\Omega^-} } w_\Omega^-[S_{n,\lambda} ,\tau_j]=\partial_{\nu_{\Omega^-} } w_\Omega^-[S_{n,\lambda} ,\tau]
\qquad\text{in}\ V^{-1,\beta} (\partial\Omega)\,.
\end{eqnarray}
Since  $\tau_j\in C^{1,\alpha}(\partial\Omega)$, it is classically known that
\[
\partial_\nu w_\Omega^+[S_{n,\lambda} ,\tau_j]=-\partial_{\nu_{\Omega^-} } w_\Omega^-[S_{n,\lambda} ,\tau_j]\qquad \text{on}\ \partial\Omega
\]
for all $j\in {\mathbb{N}}$ (cf.~\textit{e.g.} Costabel \cite[Lem.~4.1]{Co88}). Hence, the limiting relations (\ref{thm:junodedl2}) imply the validity of equality (\ref{thm:junodedl1}).\hfill  $\Box$ 

\vspace{\baselineskip}

\section{A distributional form of a   quasi-sym\-metri\-zation principle}
\label{sec:qsph}

We first prove the following classical form of a quasi-symmetrization principle that generalizes the 
Plemelj symmetrization principle, that we prove by modifying  an argument of Khavinson, Putinar and Shapiro~\cite[Lem.~2]{KhPuSh07} in a $L^p$ setting and of  Mitrea and Taylor \cite[(7.41)]{MitTa99} in a Sobolev space setting.
\begin{lemma}
 \label{lem:qlcomvw}
Let $\alpha\in ]0,1[$. Let $\Omega$ be a bounded open  subset of ${\mathbb{R}}^{n}$ of class $C^{1,\alpha}$. Let $\lambda\in {\mathbb{C}}$.  Let $S_{n,\lambda} $ be a fundamental solution 
of the operator $\Delta+\lambda$. If 
$\eta\in C^{0,\alpha}(\partial\Omega)$, then
\begin{eqnarray}\label{lem:qlcomvw1}
\lefteqn{V_\Omega[W_\Omega^t[S_{n,\lambda} ,\eta]](x)-W_\Omega[V_\Omega[S_{n,\lambda} ,\eta]](x)
}
\\ \nonumber
&&\ \ 
=\frac{1}{2}V_\Omega[\eta](x)-\frac{1}{2}V_\Omega[S_{n,\lambda} ,\eta](x)
-\lambda\int_\Omega S_n(x-y)v_\Omega^+[S_{n,\lambda} ,\eta](y)dy
\quad\forall x\in\partial\Omega\,,
\end{eqnarray}
(cf.~(\ref{eq:vosn}), (\ref{eq:wosn})).
\end{lemma}
{\bf Proof.} Let $u\equiv v_\Omega^+[S_{n,\lambda} ,\eta]$. By classical results on the acoustic  single layer potential, we have $v_\Omega^+[S_{n,\lambda} ,\eta]\in C^{1,\alpha}(\overline{\Omega})\cap C^2(\Omega)$ (cf.~\textit{e.g.}, \cite[Thm.~7.1]{DoLa17}). Then by applying the second Green Identity for the Laplace operator to the functions $v_\Omega^+[S_{n,\lambda} ,\eta]$ and $S_n(x-\cdot)$ in $\Omega$, we have
\begin{eqnarray*}
\lefteqn{
\int_\Omega v_\Omega^+[S_{n,\lambda} ,\eta](y)\Delta_yS_n(x-y)-\Delta v_\Omega^+[S_{n,\lambda} ,\eta](y)S_n(x-y)\,dy
}
\\ \nonumber
&&\qquad\qquad\qquad
=\int_{\partial\Omega}\frac{\partial}{\partial \nu_{\Omega,y}}\left(S_n(x-y)\right)V_\Omega[S_{n,\lambda} ,\eta](y)
\,d\sigma_y
\\ \nonumber
&&\qquad\qquad\qquad\quad 
-\int_{\partial\Omega}S_n(x-y)\frac{\partial}{\partial \nu_{\Omega}}v_\Omega^+[S_{n,\lambda} ,\eta](y)\,dy 
\end{eqnarray*}
for all $ x\in\Omega^-$ (cf.~\textit{e.g.}, \cite[Thm.~4.3]{DaLaMu21}). Hence, 
\begin{eqnarray*}
\lefteqn{
\int_\Omega  \lambda v_\Omega^+[S_{n,\lambda} ,\eta](y)S_n(x-y)\,dy
}
\\ \nonumber
&&\qquad\qquad\qquad   
=w_\Omega^-[V_\Omega[S_{n,\lambda} ,\eta]](x)-v_\Omega^-\left[\frac{\partial}{\partial \nu_{\Omega}}v_\Omega^+[S_{n,\lambda} ,\eta]\right](x)\qquad\forall x\in\Omega^-\,.
\end{eqnarray*}
Then the jump formulas for the harmonic single and double layer potential and the continuity of the harmonic volume potential in the whole of ${\mathbb{R}}^n$  imply that
\begin{eqnarray*}
\lefteqn{
\int_\Omega  \lambda v_\Omega^+[S_{n,\lambda} ,\eta](y)S_n(x-y)\,dy
}
\\ \nonumber
&&=
-\frac{1}{2}V_\Omega[S_{n,\lambda} ,\eta](x)+W_\Omega[V_\Omega[S_{n,\lambda} ,\eta]](x)-V_\Omega\left[-\frac{1}{2}\eta+W_\Omega^t[S_{n,\lambda} ,\eta]\right](x) 
\end{eqnarray*}
for all $x\in\partial\Omega$ and accordingly the formula of the statement holds true.

\hfill  $\Box$ 

\vspace{\baselineskip}

Next we prove a distributional form of the quasi-symmetrization principle of Lemma \ref{lem:qlcomvw}.
\begin{lemma}\label{lem:qldcomvw}
Let $\alpha\in ]0,1[$. Let $\Omega$ be a bounded open  subset of ${\mathbb{R}}^{n}$ of class $C^{1,\alpha}$. Let $\lambda\in {\mathbb{C}}$.  Let $S_{n,\lambda} $ be a fundamental solution 
of the operator $\Delta+\lambda$. If 
$\tau\in V^{-1,\alpha}(\partial\Omega)$, then 
\begin{eqnarray}\label{lem:qldcomvw1}
\lefteqn{V_\Omega[W_\Omega^t[S_{n,\lambda} ,\tau]](x)-W_\Omega[V_\Omega[S_{n,\lambda} ,\tau]](x)
}
\\ \nonumber
&&\ \ 
=\frac{1}{2}V_\Omega[\tau](x)-\frac{1}{2}V_\Omega[S_{n,\lambda} ,\tau](x)
-\lambda\int_\Omega S_n(x-y)v_\Omega^+[S_{n,\lambda} ,\tau](y)dy
\quad\forall x\in\partial\Omega\,.
\end{eqnarray}
\end{lemma}
 {\bf Proof.} By the approximation Lemma \ref{lem:apr-1a}, there exists a sequence $\{\tau_j\}_{j\in {\mathbb{N}}}$ in 
 $C^{1,\alpha}(\partial\Omega)$ as in (\ref{lem:apr-1a1}). By the quasi-symmetrization principle of Lemma \ref{lem:qlcomvw}, we know that equality (\ref{lem:qlcomvw1}) holds for $\eta=\tau_j$ for all $j\in {\mathbb{N}}$ and we now turn to show that we can take the limit as $j$ tends to $\infty$. To do so, we fix $\beta\in]0,\alpha[$ and we note that the continuity of $W_\Omega^t[S_{n,\lambda} ,\cdot]$ in $V^{-1,\beta}(\partial\Omega)$ implies that
\[
\lim_{j\to\infty}W_\Omega^t[S_{n,\lambda} ,\tau_j]=W_\Omega^t[S_{n,\lambda} ,\tau]\quad\text{in}\ V^{-1,\beta}(\partial\Omega) 
\]
(cf.~Proposition \ref{prop:cowt-1a}) and that accordingly the continuity of $V_\Omega$ from  $V^{-1,\beta}(\partial\Omega)$ to $C^{0,\beta}(\partial\Omega)$ implies that
\begin{equation}\label{lem:qldcomvw2}
\lim_{j\to\infty}V_\Omega[W_\Omega^t[S_{n,\lambda} ,\tau_j]]=V_\Omega[W_\Omega^t[S_{n,\lambda} ,\tau]]
\quad\text{in}\ C^{0,\beta}(\partial\Omega)\,,
\end{equation}
(cf.~Theorem  \ref{thm:slhco-1a}). 
 Next we note that the continuity of $V_\Omega$ from  $V^{-1,\beta}(\partial\Omega)$ to $C^{0,\beta}(\partial\Omega)$  of Theorem  \ref{thm:slhco-1a} implies that
 \[
\lim_{j\to\infty} V_\Omega[S_{n,\lambda} ,\tau_j]=V_\Omega[S_{n,\lambda} ,\tau]\quad\text{in}\ C^{0,\beta}(\partial\Omega)\,,
 \]
  and that accordingly the continuity of $W_\Omega$ from  $C^{0,\beta}(\partial\Omega)$ to itself implies that
  \begin{equation}\label{lem:qldcomvw3}
\lim_{j\to\infty}W_\Omega[V_\Omega[S_{n,\lambda} ,\tau_j]]=W_\Omega[V_\Omega[S_{n,\lambda} ,\tau]]
\quad\text{in}\ C^{0,\beta}(\partial\Omega)\,,
\end{equation}
(cf.~Theorem \ref{thm:dlay}).  Next we note that the continuity of $V_\Omega$ from  $V^{-1,\beta}(\partial\Omega)$ to $C^{0,\beta}(\partial\Omega)$ implies that
\begin{equation}\label{lem:qldcomvw3a}
\lim_{j\to\infty}V_\Omega[\tau_j]=V_\Omega[\tau]
\quad\text{in}\ C^{0,\beta}(\partial\Omega)\,,
\end{equation}
(cf.~Theorem  \ref{thm:slhco-1a}). Next we note that the continuity of $v_\Omega^+[S_{n,\lambda} ,\cdot]$ from  $V^{-1,\beta}(\partial\Omega)$ to $C^{0,\beta}(\overline{\Omega})$  of Theorem  \ref{thm:slhco-1a} implies that
 \[
\lim_{j\to\infty}v_\Omega^+[S_{n,\lambda} ,\tau_j]=v_\Omega^+[S_{n,\lambda} ,\tau]\quad\text{in}\ C^{0,\beta}(\overline{\Omega})\,,
\]
 and that accordingly the continuity of the harmonic volume potential from  $C^{0,\beta}(\overline{\Omega})$ to $C^{2,\beta}(\overline{\Omega})$ implies that
 \begin{equation}\label{lem:qldcomvw4}
\lim_{j\to\infty}\int_\Omega S_n(\cdot-y)v_\Omega^+[S_{n,\lambda} ,\tau_j](y)dy
=
\int_\Omega S_n(\cdot-y)v_\Omega^+[S_{n,\lambda} ,\tau](y)dy\quad\text{in}\ C^{2,\beta}(\overline{\Omega})\,,
\end{equation}
(cf.~Theorem \ref{thm:nwtdma} (i)).  Then the validity of equality (\ref{lem:qldcomvw1}) for $\tau_j$ for all $j\in {\mathbb{N}}$ and the limiting relations (\ref{lem:qldcomvw2})--(\ref{lem:qldcomvw4}) imply the validity of equality (\ref{lem:qldcomvw1}) for $\tau$.\hfill  $\Box$ 

\vspace{\baselineskip}

\section{A compactness result for $W_\Omega^t[S_{n,\lambda} ,\cdot]$ in $V^{-1,\alpha}(\partial\Omega)$} 
\label{sec:cowht-1a}
We first state the following embedding lemma and prove the following two regularization statements.
\begin{lemma}\label{lem:amc0b-1a}
 Let $\alpha \in]0,1[$.  Let $\Omega$ be a bounded open subset of ${\mathbb{R}}^n$ of class $C^{1,\alpha}$. Then $C^{0}(\overline{\Omega})$ is continuously embedded into $C^{-1,\alpha}(\overline{\Omega})$.
 \end{lemma}
 For a proof we refer to \cite[Lem.~21]{La24d}.
\begin{proposition}\label{prop:wh-w0b0a}
Let  $\alpha,\beta\in]0,1[$, $\beta\leq\alpha$.  Let  $\Omega$ be a  bounded open subset of ${\mathbb{R}}^{n}$ of class $C^{1,\alpha}$. Let $\lambda\in {\mathbb{C}}$.  Let $S_{n,\lambda} $ be a fundamental solution 
of the operator $\Delta+\lambda$. Then the difference
\begin{equation}\label{prop:wh-w0b0a1}
w_\Omega^+[S_{n,\lambda} ,\cdot]-w_\Omega^+[\cdot]
\end{equation}
is linear and continuous from $C^{0,\beta}(\partial\Omega)$ to $C^{0,\alpha}(\overline{\Omega})_\Delta$
\end{proposition}
{\bf Proof.} To shorten our notation, we set
\[
G[\cdot]\equiv w_\Omega^+[S_{n,\lambda} ,\cdot]-w_\Omega^+[\cdot]\,.
\]
 Next we note that
\[
\Delta (G[\mu])=-\lambda w_\Omega^+[S_{n,\lambda} ,\mu]  \ \  \text{in}\ \Omega\,,
\]
 and that accordingly, 
\begin{equation}\label{prop:wh-w0b0a2}
\left\{
\begin{array}{ll}
 \Delta (G[\mu]+{\mathcal{P}}_\Omega^+[S_n,\lambda w_\Omega^+[S_{n,\lambda} ,\mu]])=0 &\text{in}\ \Omega\,,
 \\
 G[\mu](x)+{\mathcal{P}}_\Omega^+[S_n,\lambda w_\Omega^+[S_{n,\lambda} ,\mu]](x)
 \\
 \qquad\ 
 =w_\Omega^+[S_{n,\lambda} ,\mu] (x)-w_\Omega^+[\mu](x)+{\mathcal{P}}_\Omega^+[S_n,\lambda w_\Omega^+[S_{n,\lambda} ,\mu]](x)
 &\forall x\in\partial\Omega
\end{array}
\right.
\end{equation}
for all $\mu\in C^{0,\beta}(\partial\Omega)$ (cf.~Theorem~\ref{thm:nwtdma}). By a known regularization result for the double layer potential on the boundary, we have
\[
W_\Omega [S_{n,\lambda} ,\cdot] \in {\mathcal{L}}\left(C^{0,\beta}(\partial\Omega),C^{0,\alpha}(\partial\Omega)\right)\,,
\quad W_\Omega [S_{n} ,\cdot] \in {\mathcal{L}}\left(C^{0,\beta}(\partial\Omega),C^{0,\alpha}(\partial\Omega)\right)\,,
\]
(cf.~\textit{e.g.},  \cite[Thm.~7.4]{DoLa17}) and accordingly the jump formula for the double layer potential implies that
\begin{equation}\label{prop:wh-w0b0a3}
 w_\Omega^+[S_{n,\lambda} ,\cdot]_{|\partial\Omega}(x)-w_\Omega^+[\cdot]_{|\partial\Omega}
 =W_\Omega [S_{n,\lambda} ,\cdot]-W_\Omega [S_{n} ,\cdot]\in {\mathcal{L}}\left(C^{0,\beta}(\partial\Omega),C^{0,\alpha}(\partial\Omega)\right) 
\end{equation}
(cf.~(\ref{thm:dlay1a})).  By Theorem \ref{thm:dlay} on the double layer potential, by Theorem  \ref{thm:nwtdma} on the volume potential and by the continuity of the embedding of $C^{2,\beta}(\overline{\Omega})$ into $C^{0,\alpha}(\overline{\Omega})_\Delta$ that follows by Lemma \ref{lem:amc0b-1a}, we have
\begin{equation}\label{prop:wh-w0b0a4}
{\mathcal{P}}_\Omega^+[S_n,\lambda w_\Omega^+[\cdot]]\in {\mathcal{L}}\left(C^{0,\beta}(\partial\Omega),C^{0,\alpha}(\overline{\Omega})_\Delta\right)\,.
\end{equation}
By classical properties of the Green operator, we have
\begin{equation}\label{prop:wh-w0b0a5}
{\mathcal{G}}_{\Omega,d,+}\in{\mathcal{L}}\left(C^{0,\alpha}(\partial\Omega),  C^{0,\alpha}_h(\overline{\Omega})\right)
\end{equation}
and 
\begin{eqnarray}\label{prop:wh-w0b0a6}
\lefteqn{G[\mu]=-{\mathcal{P}}_\Omega^+[S_n,\lambda w_\Omega^+[\mu]]
}
\\
\nonumber
&&
\qquad+{\mathcal{G}}_{\Omega,d,+}
\left[
w_\Omega^+[S_{n,\lambda} ,\mu]_{|\partial\Omega} -w_\Omega^+[\mu]_{|\partial\Omega} +{\mathcal{P}}_\Omega^+[S_n,\lambda w_\Omega^+[\mu]]_{|\partial\Omega}
\right]
\end{eqnarray}
for all $\mu\in C^{0,\beta}(\partial\Omega)$ (cf.~(\ref{prop:wh-w0b0a2}) and Theorem \ref{thm:idwp}). By the memberships of (\ref{prop:wh-w0b0a3})--(\ref{prop:wh-w0b0a5}) and by equality (\ref{prop:wh-w0b0a6}), we have
 \[
G[\cdot]\in {\mathcal{L}}\left(C^{0,\beta}(\partial\Omega),C^{0,\alpha}_\Delta(\overline{\Omega})\right)
\]
and thus the proof is complete.\hfill  $\Box$ 

\vspace{\baselineskip}

\begin{proposition}\label{prop:vh-v0b0a}
 Let  $\alpha,\beta\in]0,1[$, $\beta\leq\alpha$.  Let  $\Omega$ be a  bounded open subset of ${\mathbb{R}}^{n}$ of class $C^{1,\alpha}$. Let $\lambda\in {\mathbb{C}}$.  Let $S_{n,\lambda} $ be a fundamental solution 
of the operator $\Delta+\lambda$. Then the difference
\begin{equation}\label{prop:vh-v0b0a1}
v_\Omega^+[S_{n,\lambda} ,\cdot]-v_\Omega^+[\cdot]
\end{equation}
is linear and continuous from $V^{-1,\beta}(\partial\Omega)$ to $C^{0,\alpha}(\overline{\Omega})_\Delta$.
\end{proposition}
{\bf Proof.} By the definition of $V^{-1,\beta}(\partial\Omega)$, by the Lemma 13.5  of
 \cite{La24b}
 on the continuity of linear maps defined on
 $V^{-1,\beta}(\partial\Omega)$, 
 it suffices to show that if $\mu_0$, $\mu_1\in C^{0,\beta}(\partial\Omega)$, then
\begin{equation}\label{prop:vh-v0b0a2}
 v_\Omega^+[S_{n,\lambda} ,\mu_0+S_{\Omega,+}^t[\mu_1]]
-
v_\Omega^+[\mu_0+S_{\Omega,+}^t[\mu_1]]\in C^{0,\alpha}(\overline{\Omega})_\Delta\,
\end{equation}
and the map from $\left(C^{0,\beta}(\partial\Omega)\right)^2$ to $C^{0,\alpha}(\overline{\Omega})_\Delta$  that takes  $(\mu_0,\mu_1)$ to 
\[
v_\Omega^+[S_{n,\lambda} ,\mu_0+S_{\Omega,+}^t[\mu_1]]
-
v_\Omega^+[\mu_0+S_{\Omega,+}^t[\mu_1]]
\]
 is continuous. By a classical result on the single layer potential, we know that  
\begin{equation}\label{prop:vh-v0b0a3}
v_\Omega^+[S_{n,\lambda} ,\mu_0]
-
v_\Omega^+[\mu_0 ]\in {\mathcal{L}}\left(C^{0,\beta}(\partial\Omega),C^{1,\beta}(\overline{\Omega})\right)\,,
\end{equation}
(cf.~\textit{e.g.}, \cite[Thm.~7.1]{DoLa17}). Since
\[
\Delta \left[ v_\Omega^+[S_{n,\lambda} ,\mu_0 ]
-
v_\Omega^+[\mu_0]
\right]
 =-\lambda v_\Omega^+[S_{n,\lambda} ,\mu_0]  \ \  \text{in}\ \Omega\,,
\]
for all $\mu_0\in C^{0,\beta}(\partial\Omega)$, the same classical result on the single layer potential implies that
\begin{equation}\label{prop:vh-v0b0a3a}
\Delta \left[v_\Omega^+[S_{n,\lambda} ,\cdot]
-
v_\Omega^+[\cdot]\right] \in {\mathcal{L}}\left(C^{0,\beta}(\partial\Omega),C^{1,\beta}(\overline{\Omega})\right)\,.
\end{equation}
Then the memberships of (\ref{prop:vh-v0b0a3}), (\ref{prop:vh-v0b0a3a}) and the continuity of the embedding of $C^{1,\beta}(\overline{\Omega})$ into $C^{0,\alpha}(\overline{\Omega})$ and of $C^{0,\alpha}(\overline{\Omega})$ into $C^{-1,\alpha}(\overline{\Omega})$ that follows by Lemma \ref{lem:amc0b-1a} imply that 
\begin{equation}\label{prop:vh-v0b0a3b}
v_\Omega^+[S_{n,\lambda} ,\cdot]
-
v_\Omega^+[\cdot]\in {\mathcal{L}}\left(C^{0,\beta}(\partial\Omega),C^{0,\alpha}(\overline{\Omega})_\Delta\right)\,.
\end{equation}
Next we note that Theorem \ref{thm:rfdlap} implies that
\begin{eqnarray}\label{prop:vh-v0b0a4}
\lefteqn{
v_\Omega^+[S_{n,\lambda} ,S_{\Omega,+}^t[\mu_1]](x)
-
v_\Omega^+[ S_{\Omega,+}^t[\mu_1]]
}
\\ \nonumber
&&\qquad\qquad\qquad
=-{\mathcal{G}}_{d,+}[\mu_1](x)+\lambda\int_{\Omega}{\mathcal{G}}_{d,+}[\mu_1](y)S_{n,\lambda} (x-y)\,dy
\\ \nonumber
&&\qquad\qquad\qquad
+\int_{\partial\Omega}\mu_1(y)\frac{\partial}{\partial\nu_{\Omega,y}}\left(S_{n,\lambda} (x-y)\right)\,d\sigma_y
\\ \nonumber
&&\qquad\qquad\qquad
+{\mathcal{G}}_{d,+}[\mu_1](x)\\ \nonumber
&&\qquad\qquad\qquad
-\int_{\partial\Omega}\mu_1(y)\frac{\partial}{\partial\nu_{\Omega,y}}\left(S_n(x-y)\right)\,d\sigma_y
\qquad\forall x\in\Omega\,,
\end{eqnarray}
for all $\mu_1\in C^{0,\beta}(\partial\Omega)$. By a classical property of the  volume potential, we have
\[ 
{\mathcal{P}}_\Omega^+[S_{n,\lambda} ,\cdot] \in 
{\mathcal{L}}\left(C^{0,\beta}(\overline{\Omega}),C^{2,\beta}(\overline{\Omega})\right)\,,
\] 
(cf.~Theorem \ref{thm:nwtdma} (i)). Since ${\mathcal{G}}_{d,+}$ is linear and continuous from $C^{0,\beta}(\partial\Omega)$ to $C^{0,\beta}(\overline{\Omega})$ and $C^{2,\beta}(\overline{\Omega})$ is continuously embedded into 
$C^{0,\alpha}(\overline{\Omega})$ and Lemma \ref{lem:amc0b-1a} implies that $C^{2,\beta}(\overline{\Omega})$ is continuously embedded into
$C^{0,\alpha}(\overline{\Omega})_\Delta$, we conclude that
\begin{equation}\label{prop:vh-v0b0a5}
{\mathcal{P}}_\Omega^+[S_{n,\lambda} ,{\mathcal{G}}_{d,+}[\cdot]]
\in 
{\mathcal{L}}\left(C^{0,\beta}(\partial\Omega),C^{0,\alpha}(\overline{\Omega})_\Delta\right)\,.
\end{equation}
Then the regularization Proposition \ref{prop:wh-w0b0a}  implies that
\begin{equation}\label{prop:vh-v0b0a6}
w_\Omega^+[S_{n,\lambda} ,\cdot]-w_\Omega^+[S_n,\cdot]\in 
{\mathcal{L}}\left(C^{0,\beta}(\partial\Omega),C^{0,\alpha}(\overline{\Omega})_\Delta\right)\,.
\end{equation}
Hence, the formula in  (\ref{prop:vh-v0b0a4}) and the continuity properties of (\ref{prop:vh-v0b0a5})--(\ref{prop:vh-v0b0a6}) imply that
\begin{eqnarray}\label{prop:vh-v0b0a7}
v_\Omega^+[S_{n,\lambda} ,S_{\Omega,+}^t[\cdot]]
-v_\Omega^+[S_{\Omega,+}^t[\cdot]]
\in 
{\mathcal{L}}\left(C^{0,\beta}(\partial\Omega),C^{0,\alpha}(\overline{\Omega})_\Delta\right)
\,.
\end{eqnarray}
Then the continuity properties of (\ref{prop:vh-v0b0a3b}) and (\ref{prop:vh-v0b0a7}) imply the validity of (\ref{prop:vh-v0b0a2}) and of the subsequent continuity properties. Hence,  the proof is complete.\hfill  $\Box$ 

\vspace{\baselineskip}

 Next we prove the following statement that extends to the Helmholtz operator a result that is known to hold for the Laplace operator (see \cite[Thm.~19.1]{La24b}). 
\begin{theorem}\label{thm:wt-1ach}
 Let  $\alpha,\beta\in]0,1[$, $\beta\leq\alpha$.  Let  $\Omega$ be a  bounded open subset of ${\mathbb{R}}^{n}$ of class $C^{1,\alpha}$. Let $\lambda\in {\mathbb{C}}$.  Let $S_{n,\lambda} $ be a fundamental solution 
of the operator $\Delta+\lambda$. 
 Let $W_{\Omega}^t[S_{n,\lambda} ,\cdot]$ be the transpose operator to $W_{\Omega}[S_{n,\lambda} ,\cdot]$ in $C^{1,\alpha} (\partial\Omega)$. 
 Then the restriction of the operator $W_{\Omega}^t[S_{n,\lambda} ,\cdot]$ to $V^{-1,\beta}(\partial\Omega)$ is linear and continuous from  $V^{-1,\beta}(\partial\Omega)$ to $V^{-1,\alpha}(\partial\Omega)$.
\end{theorem}
{\bf Proof.} We proceed as in the proof of \cite[Thm.~19.1]{La24b}. By  \cite[Thm.~18.3]{La24b}, the operator 
 $J^{-1,\gamma}$  
from  $V^{-1,\gamma}(\partial\Omega)$ onto $C^{0,\gamma}(\partial\Omega)$ that is defined by 
\begin{equation}\label{thm:v-1an1}
J^{-1,\beta}[\tau]\equiv V_\Omega\left[\tau-\frac{\langle \tau,1\rangle }{\langle 1,1\rangle }1\right]
+
\frac{\langle \tau,1\rangle }{\langle 1,1\rangle }1 
\qquad\forall \tau\in V^{-1,\gamma}(\partial\Omega)\,.
\end{equation}
is a linear homeomorphism for all $\gamma\in]0,\alpha]$. Next, we show that $J^{-1,\beta}W_{\Omega}^t[S_{n,\lambda} ,\cdot]$ is linear and continuous from 
$V^{-1,\beta}(\partial\Omega)$ to $C^{0,\alpha}(\partial\Omega)$. To do so, we show that $J^{-1,\beta}W_{\Omega}^t[S_{n,\lambda} ,\cdot]$ equals a continuous operator from 
$V^{-1,\beta}(\partial\Omega)$ to $C^{0,\alpha}(\partial\Omega)$. Let
 $\tau\in V^{-1,\beta}(\partial\Omega)$. Then  the quasi-symmetrization of Lemma \ref{lem:qldcomvw} implies that
\begin{eqnarray}\label{thm:wt-1ach1}
\lefteqn{
 J^{-1,\beta}W_{\Omega}^t[S_{n,\lambda} ,\tau]
}
\\ \nonumber
&& 
=
V_\Omega\left[W_{\Omega}^t[S_{n,\lambda} ,\tau]-\frac{\langle W_{\Omega}^t[S_{n,\lambda} ,\tau],1\rangle }{\langle 1,1\rangle }1\right]
+
\frac{\langle W_{\Omega}^t[S_{n,\lambda} ,\tau]],1\rangle }{\langle 1,1\rangle }1
\\ \nonumber
&& 
=
V_\Omega\left[W_{\Omega}^t[S_{n,\lambda} ,\tau]\right]
-V_\Omega\left[\frac{\langle W_{\Omega}^t[S_{n,\lambda} ,\tau],1\rangle }{\langle 1,1\rangle }1\right]
+
\frac{\langle W_{\Omega}^t[S_{n,\lambda} ,\tau],1\rangle }{\langle 1,1\rangle }1
\\ \nonumber
&& 
=
W_\Omega\left[V_{\Omega}[S_{n,\lambda} ,\tau]\right]
+\frac{1}{2}V_\Omega[\tau] -\frac{1}{2}V_\Omega[S_{n,\lambda} ,\tau] 
-\lambda\int_\Omega S_n(\cdot-y)v_\Omega^+[S_{n,\lambda} ,\tau](y)dy
\\ \nonumber
&&\quad
-\frac{\langle \tau,W_{\Omega}[S_{n,\lambda} ,1]\rangle }{\langle 1,1\rangle }V_\Omega\left[1\right]
+
\frac{\langle \tau,W_{\Omega}[S_{n,\lambda} ,1]\rangle }{\langle 1,1\rangle }1
 \,.
\end{eqnarray}
By Theorem \ref{thm:slhco-1a}, $V_{\Omega}[S_{n,\lambda} ,\cdot]$ is linear and continuous from 
$V^{-1,\beta}(\partial\Omega)$ to $C^{0,\beta}(\partial\Omega)$ and  
$W_\Omega[\cdot]$ is linear and continuous from $C^{0,\beta}(\partial\Omega)$ to $C^{0,\alpha}(\partial\Omega)$  (cf. \textit{e.g.},  \cite[Thm.~4.33 (ii)]{DaLaMu21}). Hence, $W_\Omega\left[V_{\Omega}[S_{n,\lambda} ,\cdot]\right]$ is linear and continuous from $V^{-1,\beta}(\partial\Omega)$ to $C^{0,\alpha}(\partial\Omega)$. By Proposition \ref{prop:vh-v0b0a}, the difference 
$V_{\Omega}[\cdot]- V_{\Omega}[S_{n,\lambda} ,\cdot]$ is linear and continuous from 
$V^{-1,\beta}(\partial\Omega)$ to $C^{0,\alpha}(\partial\Omega)$. 

By Theorem \ref{thm:slhco-1a}, $v_{\Omega}^+[S_{n,\lambda} ,\cdot]$ is linear and continuous from the space 
$V^{-1,\beta}(\partial\Omega)$ to $C^{0,\beta}(\overline{\Omega})$ and   Theorem \ref{thm:nwtdma} implies that the harmonic volume potential is linear and continuous from $C^{0,\beta}(\overline{\Omega})$ to $C^{2,\beta}(\overline{\Omega})$, that is continuously embedded into $C^{0,\alpha}(\overline{\Omega})$. 
Hence, the map that takes $\tau$ to $\lambda\int_\Omega S_n(\cdot-y)v_\Omega^+[S_{n,\lambda} ,\tau](y)dy$
is linear and continuous from $V^{-1,\beta}(\partial\Omega)$ to $C^{0,\alpha}(\partial\Omega)$.

Since   $W_{\Omega}[S_{n,\lambda} ,1]\in C^{1,\beta}(\partial\Omega)$ (cf.~\textit{e.g.}, \cite[Thm.~7.3]{DoLa17}),   the map from $V^{-1,\beta}(\partial\Omega)$  to ${\mathbb{R}}$ that takes $\tau$ to $\frac{\langle \tau,W_{\Omega}[S_{n,\lambda} ,1]\rangle }{\langle 1,1\rangle }$ is linear and continuous. Since $V_\Omega\left[1\right]\in C^{1,\beta}(\partial\Omega)$  (cf.~\textit{e.g.}, \cite[Thm.~7.1]{DoLa17}), we conclude that the map from $V^{-1,\beta}(\partial\Omega)$  to $C^{0,\alpha}(\partial\Omega)$ that takes $\tau$ to the right hand side of (\ref{thm:wt-1ach1}) is continuous.

 Next we prove that  $W_\Omega^t[S_{n,\lambda} ,V^{-1,\beta}(\partial\Omega)]$ is contained in $V^{-1,\alpha}(\partial\Omega)$. To do so, we note that
 \begin{equation}\label{thm:wt-1ach2}
J^{-1,\beta}[\mu]=J^{-1,\alpha}[\mu]\qquad\forall \mu\in V^{-1,\alpha}(\partial\Omega)\,,
\end{equation}
that $J^{-1,\beta}$ is a bijection from $V^{-1,\beta}(\partial\Omega)$ onto $C^{0,\beta}(\partial\Omega)$, that $C^{0,\alpha}(\partial\Omega)$ is a subset of    $ C^{0,\beta}(\partial\Omega)$, that $J^{-1,\alpha}$ is a bijection from $V^{-1,\alpha}(\partial\Omega)$ onto $C^{0,\alpha}(\partial\Omega)$ and that accordingly the inclusion of $J^{-1,\beta}W_\Omega^t[S_{n,\lambda} ,V^{-1,\beta}(\partial\Omega)]$ in $C^{0,\alpha}(\partial\Omega)$ implies that  
\[
W_\Omega^t[S_{n,\lambda} , V^{-1,\beta}(\partial\Omega)]\subseteq V^{-1,\alpha}(\partial\Omega)\,.
\]
 Then equalities (\ref{thm:wt-1ach1}) and (\ref{thm:wt-1ach2}) and the above continuity of  $J^{-1,\beta}W_\Omega^t[S_{n,\lambda} ,\cdot]$ imply that $J^{-1,\alpha}W_\Omega^t[S_{n,\lambda} ,\cdot]=J^{-1,\beta}W_\Omega^t[S_{n,\lambda} ,\cdot]$ is linear and continuous from 
$V^{-1,\beta}(\partial\Omega)$ to $C^{0,\alpha}(\partial\Omega)$ and thus   Lemma 18.7 (i)  of \cite{La24b}  for linear $V^{-1,\alpha}(\partial\Omega)$-valued maps implies the continuity of $W_\Omega^t[S_{n,\lambda} ,\cdot]$ from $V^{-1,\beta}(\partial\Omega)$  to $V^{-1,\alpha}(\partial\Omega)$   and thus the proof is complete.\hfill  $\Box$ 

\vspace{\baselineskip}

By the compact embedding of $V^{-1,\alpha}(\partial\Omega)$ into $V^{-1,\beta}(\partial\Omega)$ of Theorem 18.8 of \cite{La24b}  and by Theorem \ref{thm:wt-1ach}, we conclude that the following statement holds true.
\begin{corollary}\label{corol:thm:wt-1acph}
 Let  $\alpha\in]0,1[$.  Let  $\Omega$ be a  bounded open subset of ${\mathbb{R}}^{n}$ of class $C^{1,\alpha}$. Let $\lambda\in {\mathbb{C}}$.  Let $S_{n,\lambda} $ be a fundamental solution 
of the operator $\Delta+\lambda$.  Then 
the operator $W_{\Omega}^t[S_{n,\lambda} ,\cdot]$ is compact from  $V^{-1,\alpha}(\partial\Omega)$ to itself.
\end{corollary}

\section{Appendix}
  \label{sec:appendix} 
  
  \begin{lemma}\label{lem:limprofusolh}
  Let $\lambda\in {\mathbb{C}}$.  Let $S_{n,\lambda} \in C^2({\mathbb{R}}^n\setminus\{0\})\cap L^1_{ {\mathrm{loc}} }({\mathbb{R}}^n)$ be a fundamental solution 
of the operator $\Delta+\lambda$.   Let $x\in{\mathbb{R}}^n$. Then
\begin{eqnarray}\label{lem:limprofusolh1}
\lefteqn{
\lim_{\epsilon\to 0}
\int_{\partial   {\mathbb{B}}_n(x,\epsilon)}\left[
\frac{\partial\psi}{\partial\nu_{ {\mathbb{B}}_n(x,\epsilon)} } (y) S_{n,\lambda} (x-y) -\psi (y) \frac{\partial}{\partial\nu_{{\mathbb{B}}_n(x,\epsilon),y}}\left(S_{n,\lambda} (x-y)\right)
\right]\,d\sigma
}
\\ \nonumber
&&\qquad\qquad\qquad\qquad\qquad\qquad\qquad\qquad\qquad\qquad
=-\psi(x)\qquad\forall  \psi\in C^2_c({\mathbb{R}}^n)\,.
\end{eqnarray}
\end{lemma}
{\bf Proof.} By assumption  and by the definition of derivative in the sense of distributions, we have
  \begin{eqnarray*}
\lefteqn{
  \varphi (0)= \langle(\Delta +\lambda) S_{n,\lambda}  ,\varphi  \rangle=\langle  S_{n,\lambda}  ,
 \Delta   \varphi +\lambda\varphi \rangle
}
\\ \nonumber
&&\qquad
 =\int_{{\mathbb{R}}^n} S_{n,\lambda} (y)(
 \Delta    \varphi (y)+\lambda\varphi(y)) \,dy
  \qquad\forall  \varphi \in {\mathcal{D}}({\mathbb{R}}^n)\,.
\end{eqnarray*}
 Then  the density of ${\mathcal{D}}({\mathbb{R}}^n)$ in $C^2_c({\mathbb{R}}^n)$ implies that
  \[
 \varphi(0)= \int_{{\mathbb{R}}^n} S_{n,\lambda} (y)(
 \Delta \varphi(y)+\lambda\varphi(y)) \,dy
\qquad\forall  \varphi \in C^2_c({\mathbb{R}}^n)\,.
  \]
   Next we fix $\varphi \in C^2_c({\mathbb{R}}^n)$ and we take $r\in]0,+\infty[$ such that ${\mathrm{supp}}\,\varphi\subseteq {\mathbb{B}}_n(0,r)$. Then by applying the second Green Identity in the set
\[
{\mathbb{A}}(\epsilon,r)\equiv{\mathbb{B}}_n(0,r)\setminus\overline{{\mathbb{B}}_n(0,\epsilon)} 
 \]
 to the functions
 \[
  \varphi(y)\,,\qquad  S_{n,\lambda}(y)\qquad\forall y\in 
 \overline{{\mathbb{A}}(\epsilon,r)}\,, 
 \]
 we obtain that 
 \begin{eqnarray*}
\lefteqn{
 \varphi(0)= \int_{{\mathbb{R}}^n}S_{n,\lambda} (y)
(\Delta+\lambda)\varphi(y)\,dy 
}
\\ \nonumber
&& 
=\int_{{\mathbb{B}}_n(0,\epsilon)}S_{n,\lambda} (y)
(\Delta+\lambda)\varphi(y)\,dy+\int_{{\mathbb{A}}(\epsilon,r)}S_{n,\lambda} (y)
(\Delta+\lambda)\varphi(y)\,dy
\\ \nonumber
&& 
=\int_{{\mathbb{B}}_n(0,\epsilon)}S_{n,\lambda} (y)
(\Delta+\lambda)\varphi(y)\,dy 
+
\int_{{\mathbb{A}}(\epsilon,r)}\varphi(y)(\Delta+\lambda)S_{n,\lambda} (y)\,dy 
\\ \nonumber
&&\quad
+
\int_{\partial {\mathbb{A}}(\epsilon,r)}\biggl[
\frac{\partial\varphi}{\partial\nu_{{\mathbb{A}}(\epsilon,r)}} (\eta) S_{n,\lambda} (\eta)  -\varphi  (\eta)\frac{\partial}{\partial\nu_{{\mathbb{A}}(\epsilon,r),\eta}}S_{n,\lambda} (\eta) 
\biggr]\,d\sigma_\eta
\end{eqnarray*}
for all $\epsilon\in]0,r[$. Since $\varphi$ vanishes in an open neighborhood of $\partial{\mathbb{B}}_n(0,r)$ and
\begin{eqnarray*}
&&(\Delta+\lambda)S_{n,\lambda} =0\quad\forall y\in {\mathbb{R}}^n\setminus\{0\}\,,
\\
&&\frac{\partial \varphi}{\partial\nu_{ {\mathbb{A}}(\epsilon,r)}} (\eta)=
-\frac{\partial \varphi}{\partial\nu_{  {\mathbb{B}}_n(0,\epsilon)}}  
 (\eta)\quad\forall \eta\in \partial{\mathbb{B}}_n(0,\epsilon)\,,
 \\
&&\frac{\partial S_{n,\lambda} }{\partial\nu_{ {\mathbb{A}}(\epsilon,r)}} (\eta)=
-\frac{\partial S_{n,\lambda} }{\partial\nu_{  {\mathbb{B}}_n(0,\epsilon)}}  
 (\eta)\quad\forall \eta\in \partial{\mathbb{B}}_n(0,\epsilon)\,,
 \end{eqnarray*}
we have
\begin{eqnarray*}
\lefteqn{
\varphi(0)=\int_{{\mathbb{B}}_n(0,\epsilon)}S_{n,\lambda} (\eta)
(\Delta+\lambda)\varphi(\eta)\,d\eta 
}
\\ \nonumber
&&\qquad
-
\int_{\partial {\mathbb{B}}_n(0,\epsilon) }\left[
\frac{\partial\varphi}{\partial\nu_{{\mathbb{B}}_n(0,\epsilon)}} (\eta) S_{n,\lambda} (\eta) -\varphi  (\eta)\frac{\partial}{\partial\nu_{{\mathbb{B}}_n(0,\epsilon)}}S_{n,\lambda} (\eta)  
\right] \,d\sigma_\eta
\end{eqnarray*}
for all $\epsilon\in]0,r[$. Then our assumption that $S_{n,\lambda} \in L^1_{ {\mathrm{loc}} }({\mathbb{R}}^n)$ and the Dominated Convergence Theorem imply that
\[
\lim_{\epsilon\to 0}\int_{{\mathbb{B}}_n(0,\epsilon)}S_{n,\lambda} (\eta)
(\Delta+\lambda)\varphi(\eta)\,d\eta=0\,.
\]
 Hence,  
\[
\varphi(0)=-\lim_{\epsilon\to 0}\int_{\partial {\mathbb{B}}_n(0,\epsilon) }\left[
\frac{\partial\varphi}{\partial\nu_{{\mathbb{B}}_n(0,\epsilon)}} (\eta) S_{n,\lambda} (\eta)-\varphi (\eta) \frac{\partial}{\partial\nu_{{\mathbb{B}}_n(0,\epsilon)}}S_{n,\lambda} (\eta)  
\right] \,d\sigma_\eta
 \]
Then by replacing $\varphi$ with $\psi(x-\cdot)$ and then by setting $\eta=x-y$ in the integral in the argument of the limit and by observing that
\[
\nu_{{\mathbb{B}}_n(0,\epsilon)}(x-y)
=\frac{x-y}{|x-y|}=-\frac{y-x}{|y-x|}=-\nu_{{\mathbb{B}}_n(x,\epsilon)}(y)\qquad\forall y\in \partial {\mathbb{B}}_n(x,\epsilon)
\,,
\]
 we obtain the limiting relation of the statement.\hfill  $\Box$ 

\vspace{\baselineskip}

  Next, we introduce the following variant of the third Green Identity that we need in the paper. For the convenience of the reader we include a proof.
  \begin{theorem}
\label{thm:thirdgreenh}
Let   $\lambda\in {\mathbb{C}}$. Let $S_{n,\lambda} \in C^2({\mathbb{R}}^n\setminus\{0\})\cap L^1_{ {\mathrm{loc}} }({\mathbb{R}}^n)$ be a fundamental solution 
of the operator $\Delta+\lambda$. Let $\Omega$ be a bounded open Lipschitz subset of ${\mathbb{R}}^{n}$.  If $u\in C^{1}(\overline{\Omega})\cap C^{2}(\Omega)$   and $ \Delta u\in L^{1}(\Omega)$, then the following statements hold.
\begin{enumerate}
\item[(i)]
\begin{eqnarray}
\label{thm:thirdgreenh1}
\lefteqn{
u(x)=\int_{\Omega}(\Delta u(y)+\lambda u(y))S_{n,\lambda}  (x-y)\,dy
}
\\ \nonumber
&& \ \ 
-\int_{\partial\Omega}\frac{\partial u}{\partial\nu_{\Omega}}(y)S_{n,\lambda}  (x-y)
-
u(y)\frac{\partial}{\partial\nu_{\Omega,y} }\left(S_{n,\lambda}  (x-y)\right)\,d\sigma_{y}
\quad\forall x\in\Omega\,,
\\ \label{thm:thirdgreenh1a}
\lefteqn{
0=\int_{\Omega}(\Delta u(y)+\lambda u(y))S_{n,\lambda}  (x-y)\,dy
}
\\ \nonumber
&& \ \ 
-\int_{\partial\Omega}\frac{\partial u}{\partial\nu_{\Omega}}(y)S_{n,\lambda}  (x-y)
-
u(y)\frac{\partial}{\partial\nu_{\Omega,y} }\left(S_{n,\lambda}  (x-y)\right)\,d\sigma_{y}
\quad\forall x\in\Omega^-\,.
\end{eqnarray}
\item[(ii)] If $\Delta u+\lambda u=0$ in $\Omega$, then
\begin{eqnarray}
\label{thm:thirdgreenh2}
\lefteqn{
\int_{\partial\Omega}
u(y)\frac{\partial}{\partial\nu_{\Omega,y	}}S_{n,\lambda}  (x-y)
-
\frac{\partial u}{\partial\nu_{\Omega}}(y)S_{n,\lambda}  (x-y)
\,d\sigma_{y}
}
\\ \nonumber
&&\qquad\qquad\qquad\qquad\qquad\qquad
=
\left\{
\begin{array}{ll}
u(x) & {\text{if}}\ x\in\Omega\,,
\\
0 & {\text{if}}\ x\in {\mathbb{R}}^{n}\setminus\overline{\Omega}\,.
\end{array}\right.
\end{eqnarray}
\end{enumerate}
\end{theorem}
{\bf Proof.} (i) Let $x\in\Omega$. Let $\epsilon_{0}\in]0,+\infty[$ be such that $\overline{{\mathbb{B}}_{n}(x,\epsilon_{0})}\subseteq\Omega$. We note that if $\epsilon\in]0,\epsilon_{0}]$, then $\Omega_{x,\epsilon}\equiv \Omega\setminus \overline{{\mathbb{B}}_{n}(x,\epsilon)}$ is a bounded open Lipschitz subset of ${\mathbb{R}}^{n}$. 

Since $x\notin  \overline{\Omega_{x,\epsilon}}$, we have $ S_{n,\lambda} (x-\cdot)\in C^{2}(\overline{\Omega_{x,\epsilon}})$. Moreover, $u\in C^{1}(\overline{\Omega_{x,\epsilon}})$. 
Since $S_{n,\lambda}(x-\cdot)$ is continuous and bounded in $  \overline{\Omega_{x,\epsilon}}$ and $\Delta u$ is integrable in $\Omega$, then $\Delta u (y)S_{n,\lambda} (x-y)$ is integrable in $y\in  \Omega_{x,\epsilon}$.  Since  $S_{n,\lambda} (x-y)$ satisfies the Helmholtz equation in  $y\in  \Omega_{x,\epsilon}$, then 
$  u (y)\Delta_yS_{n,\lambda} (x-y)+\lambda	u (y)  S_{n,\lambda} (x-y) =0 $ for  all $y\in  \Omega_{x,\epsilon}$. Hence,  we can apply the second Green Identity to the functions $u(y)$ and $v(y)\equiv S_{n,\lambda} (x-y)$   
  for $y$ in $\Omega_{x,\epsilon}$, and obtain
\begin{eqnarray*}
\lefteqn{
\int_{\Omega_{x,\epsilon}}(\Delta u(y)+\lambda u(y))S_{n,\lambda}  (x-y)\,dy
}
\\ \nonumber
&&\qquad
=\int_{\partial\Omega_{x,\epsilon}}\frac{\partial u}{\partial\nu_{\Omega_{x,\epsilon}}}(y)S_{n,\lambda}  (x-y)
-
u(y)\frac{\partial}{\partial\nu_{\Omega_{x,\epsilon}, y		}}\left(S_{n,\lambda}  (x-y)\right)\,d\sigma_{y}\,,
\end{eqnarray*}
(cf.~\textit{e.g.},  \cite[Thm.~4.3]{DaLaMu21}).  Next we note that
\begin{eqnarray*}
\partial\Omega_{x,\epsilon}&=&\partial\Omega\cup\partial  {\mathbb{B}}_n(x,\epsilon)\,,
\\
\nu_{\Omega_{x,\epsilon}}=\nu_{\Omega}\qquad{\mathrm{on}}\ \partial\Omega
&\qquad&
\nu_{\Omega_{x,\epsilon}}=-\nu_{{\mathbb{B}}_n(x,\epsilon)}\qquad{\mathrm{on}}\ \partial {\mathbb{B}}_n(x,\epsilon)\,,
\end{eqnarray*}
and that accordingly
\[
\frac{\partial u}{ \partial	 \nu_{ \Omega_{x,\epsilon}}} (y)=\frac{\partial u}{ \partial	 \nu_{ \Omega}} (y)
\,,\quad
\frac{\partial }{\partial	\nu_{ \Omega_{x,\epsilon},y}}\left( S_{n,\lambda} (x-y)\right)=\frac{\partial }{
\partial	\nu_{ \Omega,y}}\left( S_{n,\lambda} (x-y)\right)
\]
for all $y\in\partial\Omega$ and
\begin{eqnarray*}
\frac{\partial u}{\partial\nu_{ \Omega_{x,\epsilon}}} (y)=
-\frac{\partial u}{\partial\nu_{  {\mathbb{B}}_n(x,\epsilon)}}  (y)\,,\ \ 
\frac{\partial  }{\partial\nu_{ \Omega_{x,\epsilon},y}}\left( S_{n,\lambda} (x-y)\right)=
-\frac{\partial  }{\partial\nu_{  {\mathbb{B}}_n(x,\epsilon),y}}  
\left( S_{n,\lambda} (x-y)\right)
\end{eqnarray*}
 for all  	$ y\in \partial{\mathbb{B}}_n(x,\epsilon)$. Hence, 
\begin{eqnarray}\label{thm:thirdgreenh3}
\lefteqn{
\int_{\Omega_{x,\epsilon}}(\Delta u(y)+\lambda u(y))S_{n,\lambda}  (x-y)\,dy
}
\\ \nonumber
&& 
=\int_{\partial\Omega}\frac{\partial u}{\partial\nu_{\Omega }}(y)S_{n,\lambda}  (x-y)
-
u(y)\frac{\partial}{\partial\nu_{\Omega, y	  }}\left(S_{n,\lambda}  (x-y)\right)\,d\sigma_{y}
\\ \nonumber
&&\quad
-\int_{\partial{\mathbb{B}}_n(x,\epsilon)}\frac{\partial u}{\partial\nu_{{\mathbb{B}}_n(x,\epsilon)}} (y)S_{n,\lambda}  (x-y)
-
u(y)\frac{\partial}{\partial\nu_{{\mathbb{B}}_n(x,\epsilon),	y	} }\left(S_{n,\lambda}  (x-y)\right)\,d\sigma_{y}\,.
\end{eqnarray}
Since  $u$,	 	$\Delta u\in C^{0}(\overline{  {\mathbb{B}}_{n}(x,\epsilon_{0}) })$ and $
S_{n,\lambda} (x-\cdot)\in L^{1}( {\mathbb{B}}_{n}(x,\epsilon_{0}))$,  the function  
\[
S_{n,\lambda} (x-y)(\Delta u(y)+\lambda u(y))
\]
 is integrable in 
$y\in  {\mathbb{B}}_{n}(x,\epsilon_{0})$. On the other hand, 
$S_{n,\lambda} (x-\cdot)$ is continuous and bounded in $  \overline{\Omega_{x,
\epsilon}}$ and $(\Delta u +\lambda u )$ is integrable in $\Omega$, and accordingly,  
\[
(\Delta u(y)+\lambda u(y))S_{n,\lambda} (x-y)
\]
  is integrable in 
$\Omega\setminus {\mathbb{B}}_{n}(x,\epsilon_{0})$. Hence, the function  $ \chi_{\Omega} (y)S_{n,\lambda}(x-y)(\Delta u(y)+\lambda u(y))$ is integrable in $y\in\Omega$ and the Dominated Convergence Theorem implies that
\begin{eqnarray}\label{thm:thirdgreenh4}
\lefteqn{\lim_{\epsilon\to 0}\int_{\Omega } \chi_{\Omega_{x,\epsilon}}(y)S_{n,\lambda}(x-y)(\Delta u(y)+\lambda u(y))\,dy
}
\\ \nonumber
&&\qquad\qquad\qquad\qquad
=
\int_{\Omega } \chi_{\Omega }(y)S_{n,\lambda}(x-y)(\Delta u(y)+\lambda u(y))\,dy\,.
\end{eqnarray}
Next we take $\varphi_x\in {\mathcal{D}}(\Omega)$ such that $\varphi_x$ equals $1$ in an open neighborhood of $ {\mathbb{B}}_{n}(x,\epsilon_{0})$ and we set
\[
u_1\equiv \varphi_x u\,,\qquad u_2\equiv (1-\varphi_x) u
\]
Since $u_1\in C^2_c({\mathbb{R}}^n)$, then Lemma \ref{lem:limprofusolh} implies that
the integral on $\partial{\mathbb{B}}_{n}(x,\epsilon)$ in the right hand side of (\ref{thm:thirdgreenh3}) with $u$ replaced by $u_1$ tends to $-u_1(x)=-u(x)$ as $\epsilon $ tends to $0$. Moreover, $u_2$ equals zero in an open neighborhood of $\partial{\mathbb{B}}_{n}(x,\epsilon)$  and thus the integral on $\partial{\mathbb{B}}_{n}(x,\epsilon)$ in the right hand side of (\ref{thm:thirdgreenh3}) with $u$ replaced by $u_2$ is equals to zero for all $\epsilon\in ]0,\epsilon_0[$. Then by taking the limit as $\epsilon$ tends to zero in (\ref{thm:thirdgreenh3})
and by invoking (\ref{thm:thirdgreenh4}), we obtain equality (\ref{thm:thirdgreenh1}).

  If instead 	$x\in\Omega^-$, then $v(\cdot)\equiv S_{n,\lambda} (x-\cdot)\in C^{2}(\overline{\Omega})$ and thus the second Green Identity implies that
\[
\int_\Omega (\Delta u +\lambda u)v -u(\Delta v+\lambda v)\,dx=\int_{\partial\Omega}\frac{\partial u}{\partial\nu}v-
\frac{\partial v}{\partial\nu}u\,d\sigma\,,
\]
and accordingly equality (\ref{thm:thirdgreenh1a}) holds true (cf.~\textit{e.g.}, \cite[Thm.~4.3]{DaLaMu21}). Statement (ii) is an immediate consequence of statement (i).\hfill  $\Box$ 

\vspace{\baselineskip}

\vspace{\baselineskip}

  \noindent
{\bf Statements and Declarations}\\

 \noindent
{\bf Competing interests:} This paper does not have any  conflict of interest or competing interest.

 \noindent
{\bf Acknowledgement.}  The author  acknowledges  the support of GNAMPA-INdAM    and   of the Project funded by the European Union – Next Generation EU under the National Recovery and Resilience Plan (NRRP), Mission 4 Component 2 Investment 1.1 - Call for tender PRIN 2022 No. 104 of February, 2 2022 of Italian Ministry of University and Research; Project code: 2022SENJZ3 (subject area: PE - Physical Sciences and Engineering) ``Perturbation problems and asymptotics for elliptic differential equations: variational and potential theoretic methods''.

\end{document}